\documentclass{amsproc}
\usepackage{amssymb}
\usepackage{amsmath}
\usepackage{amsfonts}
\usepackage[all]{xy}
\usepackage{amscd}

\theoremstyle{plain}

\newtheorem{proposition}{Proposition}

\newtheorem{theorem}{Theorem}
\numberwithin{equation}{section}
\newcommand{\squarecell}[9]{
\xymatrix@R=1pc@C=5pc{
#1 \ar@{->}[r]|{#2}="2" \ar[d]_{#4} & #3 \ar[d]^{#6} \\
#7 \ar@{->}[r]|{#8}="8" & #9  \\
\ar@{=>}^{#5} "2";"8"}}
\input{tcilatex}

\begin{document}
\title[The Category of Pseudo-Categories]{Pseudo-Categories}
\author{N. Martins-Ferreira}
\address{ESTG-IPLeiria}
\email{nelsonmf@estg.ipleiria.pt}
\urladdr{http://www.estg.ipleiria.pt/\symbol{126}nelsonmf}
\thanks{The author thanks to Professor G. Janelidze for reading and
suggesting useful changes to clarify the subject, and also to Professors R.
Brown and T. Porter for the enlightning discutions during the seminars
presented at Bangor.}
\date{April 9, 2006; 22:58 (After reviewing)}
\subjclass[2000]{Primary 18D05, 18D15; Secondary 57D99.}
\keywords{Bicategory, double category, weak category, pseudo-category,
pseudo double category, pseudo-functor, pseudo-natural transformation,
pseudo-modification, weakly cartesian closed.}
\dedicatory{}
\thanks{This paper is in final form.}

\begin{abstract}
We provide a complete description of the category of pseudo-categories
(including pseudo-functors, natural and pseudo-natural transformations and
pseudo modifications). A pseudo-category is a non strict version of an
internal category. It was called a weak category and weak double category in
some earlier papers. When internal to Cat it is at the same time a
generalization of a bicategory and a double category. The category of
pseudo-categories is a kind of \textquotedblleft
tetracategory\textquotedblright\ and it turns out to be cartesian closed in
a suitable sense.
\end{abstract}

\maketitle

\section{Introduction}

The notion of pseudo-category\footnote{%
In the previous work \cite{MF2} the word "weak" was used with the same
meaning. We claim that "pseudo" is more apropriate because it is the
intermediate term between precategory and internal category. Also it agrees
with the notion of pseudo-functor, already well established.} considered in
this paper is closely related and essentially is a special case of several
higher categorical structures studied for example by Grandis and Par\'{e} 
\cite{Pare}, Leinster \cite{Leinster}, Street \cite{Street},\cite{Street2},
among several others. We have arrived to the present definition of
pseudo-category (which some authors would probably call a pseudo double
category) while describing internal bicategories in Ab \cite{MF1}. We even
found it easier, for our particular purposes, to work with pseudo-categories
than to work with bicategories. Defining a pseudo-category we begin with a
2-category, take the definition of \ an internal category there, and replace
the equalities in the associativity and identity axioms by the existence of
suitable isomorphisms which then have to satisfy some coherence conditions.
That is, let \textbf{C}\ be a 2-category, a pseudo-category in (internal to) 
\textbf{C} is a system%
\begin{equation*}
\left( C_{0},C_{1},d,c,e,m,\alpha ,\lambda ,\rho \right) 
\end{equation*}%
where $C_{0},C_{1}$ are objects of $\mathbf{C}$,%
\begin{equation*}
d,c:C_{1}\longrightarrow C_{0}\ ,\ e:C_{0}\longrightarrow C_{1}\ ,\
m:C_{1}\times _{C_{0}}C_{1}\longrightarrow C_{1}
\end{equation*}%
are morphisms of \textbf{C} , with $C_{1}\times _{C_{0}}C_{1}$ the object in
the pullback diagram%
\begin{equation*}
\begin{tabular}{ll}
$%
\begin{CD}%
\newline
C_{1}\times _{C_{0}}C_{1}%
@>%
\pi _{2}%
>%
>%
C_{1}%
\\ %
\newline
@V%
\pi _{1}%
V%
V%
@V%
V%
c%
V%
\\ %
\newline
C_{1}%
@>%
d%
>%
>%
C_{0}%
&%
&%
\\ %
\newline
\end{CD}%
$ & ;%
\end{tabular}%
\end{equation*}%
\begin{eqnarray*}
\alpha  &:&m\left( 1_{C_{1}}\times _{C_{0}}m\right) \longrightarrow m\left(
m\times _{C_{0}}1_{C_{1}}\right) \ , \\
\lambda  &:&m\left\langle ec,1_{C_{1}}\right\rangle \longrightarrow
1_{C_{1}}\ ,\ \rho :m\left\langle 1_{C_{1}},ed\right\rangle \longrightarrow
1_{C_{1}},
\end{eqnarray*}%
are 2-cells of \textbf{C} (which are isomorphisms), the following conditions
are satisfied%
\begin{equation}
de=1_{c_{0}}=ce,  \label{RGraph}
\end{equation}%
\begin{equation}
dm=d\pi _{2}\;,\;cm=c\pi _{1},  \label{dom cod for composition}
\end{equation}%
\begin{eqnarray}
d\circ \lambda  &=&1_{d}=d\circ \rho ,  \label{dom cod for lambda rho} \\
c\circ \lambda  &=&1_{c}=c\circ \rho ,  \notag
\end{eqnarray}%
\begin{equation}
d\circ \alpha =1_{d\pi _{3}}\ ,\ c\circ \alpha =1_{c\pi _{1}},
\label{dom cod for alpha}
\end{equation}%
\begin{equation}
\lambda \circ e=\rho \circ e,  \label{id for lambda rho}
\end{equation}%
and the following diagrams commute%
\begin{equation}
\xymatrix{ & \bullet  \ar[rr]^{ m \circ (1_{C_1}\times_{C_0}\alpha)}
             \ar[ldd] _{ \alpha \circ ( 1_{C_1} \times_{C_0} 1_{C_1} \times_{C_0} m ) }
 & & \bullet \ar[rdd]^{\alpha \circ  ( 1_{C_1} \times_{C_0} m \times_{C_0} 1_{C_1} ) }  & \\
& & & & \\
\bullet   \ar[rrdd]_{\alpha \circ ( m \times_{C_0} 1_{C_1} \times_{C_0} 1_{C_1} ) }
  &  &  &  &  \bullet  \ar[ddll]^{m \circ (\alpha\times_{C_0} 1_{C_1})}  \\
& & & & \\
 & & \bullet &  &}
\label{pentagon coherence}
\end{equation}%
\begin{equation}
\xymatrix@=4pc{
 \bullet   \ar[rr]^{\alpha \circ (1_{C_1}\times_{C_0} <ec , 1_{C_1}>)}
            \ar[rd]_{m \circ (1_{C_1} \times_{C_0} \lambda)}
&& \bullet   \ar[ld]^{m \circ (\rho \times_{C_0} 1_{C_1})} \\
& \bullet}
\label{triangle coherence}
\end{equation}

Examples:

\begin{enumerate}
\item When \textbf{C}=Set with the discrete 2-category structure (only
identity 2-cells) one obtains the definition of an ordinary category since $%
\alpha ,\lambda ,\rho $ are all identities;

\item When \textbf{C}=Set with the codiscrete 2-category structure (exactly
one 2-cell for each pair of morphisms) one obtain the definition of a
precategory since $\alpha ,\lambda ,\rho $ always exist and the coherence
conditions are trivially satisfied;\newline
(This result applies equally to any category)

\item When \textbf{C}=Grp considered as a 2-category: every group is a (one
object) category and the inclusion functor%
\begin{equation*}
\text{Grp}\longrightarrow \text{Cat}
\end{equation*}%
induces a 2-category structure in Grp, where a 2-cell%
\begin{equation*}
\tau :f\longrightarrow g\ \ ,\ (f,g:A\longrightarrow B\ \text{group
homomorphisms})
\end{equation*}%
is an element $\tau \in B$, such that for every $x\in A$,%
\begin{equation*}
g\left( x\right) =\tau f\left( x\right) \tau ^{-1}.
\end{equation*}%
With \ this setting, a pseudo-category in Grp is described (see \cite{MF3})
by a group homomorphism 
\begin{equation*}
\partial :X\longrightarrow B,
\end{equation*}%
an arbitrary element%
\begin{equation*}
\delta \in \ker \partial
\end{equation*}%
and an action of $B$ in $X$ (denoted by $b\cdot x$ for $b\in B$ and $x\in X$%
) satisfying%
\begin{eqnarray*}
\partial \left( b\cdot x\right) &=&b\partial \left( x\right) b^{-1} \\
\partial \left( x\right) \cdot x^{\prime } &=&x+x^{\prime }-x
\end{eqnarray*}%
for every $b\in B,\ x,x^{\prime }\in X$. Note that the difference to a
crossed module (description of an internal category in Grp) is that \ in a
crossed module the element $\delta =1.$\newline
The pseudo-category so obtained is as follows: objects are the elements of $%
B $, arrows are pairs $\left( x,b\right) :b\longrightarrow \partial x+b$ and
the composition of $\left( x^{\prime },\partial x+b\right) :\partial
x+b\longrightarrow \partial x^{\prime }+\partial x+b$ with $\left(
x,b\right) :b\longrightarrow \partial x+b$ is the pair $\left( x^{\prime
}+x-\delta +b\cdot \delta ,b\right) :b\longrightarrow \partial x^{\prime
}+\partial x+b.$ The isomorphism between $\left( 0,\partial x+b\right) \circ
\left( x,b\right) =\left( x,b\right) \circ \left( 0,b\right) $ and $\left(
x,b\right) $ is the element $\left( \delta ,0\right) \in X\rtimes B$ .
Associativity is satisfied, since $\left( x^{\prime \prime },\partial
x^{\prime }+\partial x+b\right) \circ \left( \left( x^{\prime },\partial
x+b\right) \circ \left( x,b\right) \right) =\left( \left( x^{\prime \prime
},\partial x^{\prime }+\partial x+b\right) \circ \left( x^{\prime },\partial
x+b\right) \right) \circ \left( x,b\right) .$

\item When \textbf{C}=Mor(Ab) the 2-category of morphisms of abelian groups,
the above definition gives a structure which is completely determined by a
commutative square%
\begin{equation*}
\begin{CD}%
\newline
A_{1}%
@>%
\partial 
>%
>%
A_{0}%
\\ %
\newline
@V%
k_{1}%
V%
V%
@V%
V%
k_{0}%
V%
\\ %
\newline
B_{1}%
@>%
\partial ^{\prime }%
>%
>%
B_{0}%
\\ %
\newline
\end{CD}%
\end{equation*}%
together with three morphisms%
\begin{eqnarray*}
\lambda ,\rho &:&A_{0}\longrightarrow A_{1}, \\
\eta &:&B_{0}\longrightarrow A_{1},
\end{eqnarray*}%
satisfying conditions%
\begin{eqnarray*}
k_{1}\lambda &=&0=k_{1}\rho , \\
k_{1}\eta &=&0,
\end{eqnarray*}%
and it may be viewed as a structure with objects, vertical arrows,
horizontal arrows and squares, in the following way (see \cite{MF2}, p. 409,
for more details)%
\begin{equation*}
\begin{tabular}{rcl}
$b$ & $^{\underrightarrow{\ \ \ \ \ \ \left( b,x\right) \ \ \ \ \ \ }}$ & $%
b+k_{0}\left( x\right) $ \\ 
$\left. \binom{b}{d}\right\downarrow $ & $\left( 
\begin{array}{cc}
b & x \\ 
d & y%
\end{array}%
\right) $ & $\left\downarrow \binom{b+k_{0}\left( x\right) }{d+k_{1}\left(
y\right) }\right. $ \\ 
$b+\partial ^{\prime }\left( d\right) $ & $_{\overrightarrow{\left(
b+\partial ^{\prime }\left( d\right) ,x+\partial \left( y\right) \right) }}$
& $\ast $%
\end{tabular}%
,
\end{equation*}%
where $\ast $ stands for $b+\partial ^{\prime }\left( d\right) +k_{0}\left(
x+\partial \left( y\right) \right) =b+k_{0}\left( x\right) +\partial
^{\prime }\left( d+k_{1}\left( y\right) \right) .$

\item When \textbf{C}=Top (with homotopy classes as 2-cells) we find the
following particular example. Let $X$ be a space and consider the following
diagram%
\begin{equation*}
X^{I}\times _{X}X^{I}\overset{m}{\longrightarrow }X^{I}\underset{c}{\overset{%
d}{\underrightarrow{\overrightarrow{\overset{e}{\longleftarrow }}}}}X
\end{equation*}%
where $X^{I}$ is equipped with the compact open topology and $X^{I}\times
_{X}X^{I}$ with the product topology ($I$ is the unit interval), with 
\begin{equation*}
X^{I}\times _{X}X^{I}=\left\{ \left\langle g,f\right\rangle \mid f\left(
0\right) =g\left( 1\right) \right\}
\end{equation*}%
and $d,e,c,m$ defined as follows%
\begin{eqnarray*}
d\left( f\right) &=&f\left( 0\right) \\
c\left( f\right) &=&f\left( 1\right) \\
e_{x}\left( t\right) &=&x \\
m\left( f,g\right) &=&\left\{ 
\begin{array}{c}
g\left( 2t\right) ,t<\frac{1}{2} \\ 
f\left( 2t-1\right) ,t\geq \frac{1}{2}%
\end{array}%
\right.
\end{eqnarray*}%
with \ $f,g:I\longrightarrow X$ (continuous maps)\ and $x\in X$. The
homotopies $\alpha ,\lambda ,\rho $ are the usual ones.

\item When \textbf{C}=Cat the objects $C_{0}$ and $C_{1}$ are (small)
categories, and the morphisms $d,c,e,m$ are functors. We denote the objects
of $C_{0}$ by the first capital letters in the alphabet (possible with
primes) $A,A^{\prime },B,B^{\prime },...$ and the morphisms by first small
letters in the alphabet $a:A\longrightarrow A^{\prime },b:B\longrightarrow
B^{\prime },...$ . We will denote the objects of $C_{1}$ by small letters as 
$f,f^{\prime },g,g^{\prime },...$ and the morphisms by small greek letters
as $\varphi :f\longrightarrow f^{\prime },\gamma :g\longrightarrow g^{\prime
},...$ . We will also consider that the functors $d$ and $c$ are defined as
follows%
\begin{equation*}
\begin{tabular}{clc}
$C_{1}$ &  & $C_{0}$ \\ 
& $%
\begin{array}{c}
\\ 
d\nearrow%
\end{array}%
$ & 
\begin{tabular}{l}
$a:A\longrightarrow A^{\prime }$ \\ 
\end{tabular}
\\ 
$\varphi :f\longrightarrow f^{\prime }$ &  &  \\ 
& $c\searrow $ & 
\begin{tabular}{l}
\\ 
$b:B\longrightarrow B^{\prime }$%
\end{tabular}%
\end{tabular}%
\end{equation*}%
hence, the objects of $C_{1}$ are arrows $f:A\longrightarrow B,f^{\prime
}:A^{\prime }\longrightarrow B^{\prime },$ that we will always represent
using inplace notation as $%
\xymatrix{A \ar[r]|{f} & B \\}%
,%
\xymatrix{A^{\prime} \ar[r]|{f^{\prime}} & B^{\prime} \\}%
$ to distinguish from the morphisms of $C_{0}$,\ and thus the morphisms of $%
C_{1}$ are of the form%
\begin{equation*}
\squarecell{A}   {f}    {B}
                  {a}   {\varphi}    {b}
                  {A^\prime}   {f^\prime}    {B^\prime}%
\begin{array}[t]{c}
\\ 
.%
\end{array}%
\end{equation*}%
The functor $e$ sends $a:A\longrightarrow A^{\prime }$ to 
\begin{equation*}
\xymatrix@R=1pc@C=5pc{
A \ar[r]|{id_A}="f" \ar[d]_{a} & A \ar[d]^{a} \\
A^{\prime} \ar[r]|{id_{A^{\prime}}}="k" & A^{\prime}  \\
\ar@{=>}^{id_a} "f";"k"}%
\begin{array}[t]{c}
\\ 
,%
\end{array}%
\end{equation*}%
while the functor $m$ sends $\left\langle \gamma ,\varphi \right\rangle $ to 
$\gamma \otimes \varphi $ as displayed in the diagram below%
\begin{equation*}
\xymatrix@R=1pc@C=5pc{
A \ar[r]|{f}="f" \ar[d]_{a} & B \ar[d]^{b} \ar[r]|{g}="g" & C \ar[d]^{c} \\
A^{\prime} \ar[r]|{f^{\prime}}="fp" & B^{\prime}  \ar[r]|{g^\prime}="gp" & C^{\prime}\\
\ar@{=>}^{\varphi} "f";"fp"
\ar@{=>}^{\gamma} "g";"gp"}%
\begin{array}[t]{c}
\\ 
\longmapsto%
\end{array}%
\squarecell{A}   {g \otimes f}    {C}
                  {a}   {\gamma \otimes \varphi}    {c}
                  {A^\prime}   {g^\prime \otimes f^\prime}    {C^\prime}%
\begin{array}[t]{c}
\\ 
.%
\end{array}%
\end{equation*}%
Each component of $\alpha $ is of the form%
\begin{equation*}
\squarecell{A}   {h \otimes (g \otimes f)}    {D}
                  {1_A}   {\alpha_{h,g,f}}    {1_D}
                  {A}   {(h \otimes g )\otimes f}    {D}%
\begin{array}[t]{c}
\\ 
,%
\end{array}%
\end{equation*}%
while the components of $\lambda $ and $\rho $ are given by%
\begin{equation*}
\squarecell{A}   {id_B \otimes f}    {B}
                  {1_A}   {\lambda_f}    {1_B}
                  {A}   {f}    {B}%
\begin{array}[t]{c}
\\ 
\ \ \ \ ,\ \ \ \ \ 
\end{array}%
\squarecell{A}   {f \otimes id_A}    {B}
                  {1_A}   {\rho_f}    {1_B}
                  {A}   {f}    {B}%
\begin{array}[t]{c}
\\ 
.%
\end{array}%
\end{equation*}%
Thus, a description of pseudo-category in Cat is as follows.
\end{enumerate}

A pseudo-category in Cat is a structure with\newline
- objects: $A,A^{\prime },A^{\prime \prime },B,B^{\prime },...$\newline
- morphisms: $a:A\longrightarrow A^{\prime },a^{\prime }:A^{\prime
}\longrightarrow A^{\prime \prime },b:B\longrightarrow B^{\prime },...$%
\newline
- \textit{pseudo-morphisms}: $%
\xymatrix{A \ar[r]|{f} & B \\}%
,%
\xymatrix{A^{\prime} \ar[r]|{f^{\prime}} & B^{\prime} \\}%
,%
\xymatrix{B \ar[r]|{g} & C \\}%
,...$ \newline
- and cells:%
\begin{equation*}
\squarecell{A}   {f}    {B}
                  {a}   {\varphi}    {b}
                  {A^\prime}   {f^\prime}    {B^\prime}%
\begin{array}[t]{c}
\\ 
\ ,\ 
\end{array}%
\squarecell{A^\prime}   {f^\prime}    {B^\prime}
                  {a^\prime}   {\varphi^\prime}    {b^\prime}
                  {A^{\prime\prime}}   {f^{\prime\prime}}    {B^{\prime\prime}}%
\begin{array}[t]{c}
\\ 
\ ,\ 
\end{array}%
\squarecell{B}   {g}    {C}
                  {b}   {\gamma}    {c}
                  {B^\prime}   {g^\prime}    {C^\prime}%
\begin{array}[t]{c}
\\ 
,...%
\end{array}%
\end{equation*}%
where objects and morphisms form a category 
\begin{eqnarray*}
a^{\prime \prime }\left( a^{\prime }a\right) &=&\left( a^{\prime \prime
}a^{\prime }\right) a, \\
1_{A^{\prime }}a &=&a1_{A}\text{;}
\end{eqnarray*}%
pseudo-morphisms and cells also form a category%
\begin{eqnarray*}
\varphi ^{\prime \prime }\left( \varphi ^{\prime }\varphi \right) &=&\left(
\varphi ^{\prime \prime }\varphi ^{\prime }\right) \varphi , \\
1_{f^{\prime }}\varphi &=&\varphi 1_{f}\text{,}
\end{eqnarray*}%
with $1_{f}$ being the cell%
\begin{equation*}
\squarecell{A}   {f}    {B}
                  {1_A}   {1_f}    {1_B}
                {A}   {f}    {B}%
\begin{array}[t]{c}
\\ 
;%
\end{array}%
\end{equation*}%
for each pair of \emph{pseudo-composable} cells $\gamma ,\varphi ,$ there is
a pseudo-composition $\gamma \otimes \varphi $ 
\begin{equation*}
\squarecell{A}   {g \otimes f}    {C}
                  {a}   {\gamma \otimes \varphi}    {c}
                  {A^\prime}   {g^\prime \otimes f^\prime}    {C^\prime}%
\begin{array}[t]{c}
\\ 
;%
\end{array}%
\end{equation*}%
satisfying%
\begin{eqnarray}
\left( \gamma ^{\prime }\gamma \right) \otimes \left( \varphi ^{\prime
}\varphi \right) &=&\left( \gamma ^{\prime }\otimes \varphi ^{\prime
}\right) \left( \gamma \otimes \varphi \right) ,  \label{interchange law} \\
1_{g\otimes f} &=&1_{g}\otimes 1_{f};  \notag
\end{eqnarray}%
for each morphism $a:A\longrightarrow A^{\prime }$, there is a
pseudo-identity $id_{a}$%
\begin{equation*}
\xymatrix@R=1pc@C=5pc{
A \ar[r]|{id_A}="f" \ar[d]_{a} & A \ar[d]^{a} \\
A^{\prime} \ar[r]|{id_{A^{\prime}}}="k" & A^{\prime}  \\
\ar@{=>}^{id_a} "f";"k"}%
\end{equation*}%
satisfying%
\begin{eqnarray*}
id_{1_{A}} &=&1_{id_{A}} \\
id_{a^{\prime }a} &=&id_{a^{\prime }}id_{a};
\end{eqnarray*}%
there is a special cell $\alpha _{h,g,f}$ for each triple of composable
pseudo-morphisms $h,g,f$%
\begin{equation*}
\squarecell{A}   {h \otimes (g \otimes f)}    {D}
                  {1_A}   {\alpha_{h,g,f}}    {1_D}
                  {A}   {(h \otimes g )\otimes f}    {D}%
\begin{array}[t]{c}
\\ 
,%
\end{array}%
\end{equation*}%
natural in each component, i.e., the following diagram of cells%
\begin{equation*}
\begin{tabular}{ll}
$%
\begin{CD}%
\newline
h\otimes \left( g\otimes f\right) 
@>%
\alpha _{h,g,f}%
>%
>%
\left( h\otimes g\right) \otimes f%
\\ %
\newline
@V%
\eta \otimes \left( \gamma \otimes \varphi \right) 
V%
V%
@V%
V%
\left( \eta \otimes \gamma \right) \otimes \varphi 
V%
\\ %
\newline
h^{\prime }\otimes \left( g^{\prime }\otimes f^{\prime }\right) 
@>%
\alpha _{h^{\prime },g^{\prime },f^{\prime }}%
>%
>%
\left( h^{\prime }\otimes g^{\prime }\right) \otimes f^{\prime }%
&%
&%
\\ %
\newline
\end{CD}%
$ & 
\end{tabular}%
\end{equation*}%
commutes for every triple of pseudo-composable cells $\varphi ,\gamma ,\eta $%
\begin{equation*}
\xymatrix@R=1pc@C=5pc{
A \ar[r]|{f}="f" \ar[d]_{a} & B \ar[d]^{b} \ar[r]|{g}="g" & C  \ar[r]|{h}="h" \ar[d]^{c} & D \ar[d]^d\\
A^{\prime} \ar[r]|{f^{\prime}}="fp" & B^{\prime}  \ar[r]|{g^\prime}="gp" & C^{\prime} \ar[r]|{h^\prime}="hp"  &  D^{\prime}  \\
\ar@{=>}^{\varphi} "f";"fp"
\ar@{=>}^{\gamma} "g";"gp"
\ar@{=>}^{\eta} "h";"hp"}%
\begin{array}[t]{c}
\\ 
;%
\end{array}%
\end{equation*}%
to each pseudo-morphism $f:A\longrightarrow B$ there are two special cells%
\begin{equation*}
\squarecell{A}   {id_B \otimes f}    {B}
                  {1_A}   {\lambda_f}    {1_B}
                  {A}   {f}    {B}%
\begin{array}[t]{c}
\\ 
\ \ \ \ ,\ \ \ \ \ 
\end{array}%
\squarecell{A}   {f \otimes id_A}    {B}
                  {1_A}   {\rho_f}    {1_B}
                  {A}   {f}    {B}%
\begin{array}[t]{c}
\\ 
,%
\end{array}%
\end{equation*}%
natural in $f$, that is, to every cell $\varphi $ as above, the following
diagrams of cells commute%
\begin{equation*}
\begin{tabular}{llll}
$%
\begin{CD}%
\newline
id_{B}\otimes f%
@>%
\lambda _{f}%
>%
>%
f%
\\ %
\newline
@V%
id_{1_{B}}\otimes \varphi 
V%
V%
@V%
V%
\varphi 
V%
\\ %
\newline
id_{B^{\prime }}\otimes f^{\prime }%
@>%
\lambda _{f^{\prime }}%
>%
>%
f^{\prime }%
&%
&%
\\ %
\newline
\end{CD}%
$ & \ \ \ \ ,\ \ \ \  & $%
\begin{CD}%
\newline
f\otimes id_{A}%
@>%
\rho _{f}%
>%
>%
f%
\\ %
\newline
@V%
\varphi \otimes id_{1_{A}}%
V%
V%
@V%
V%
\varphi 
V%
\\ %
\newline
f^{\prime }\otimes id_{B^{\prime }}%
@>%
\lambda _{f^{\prime }}%
>%
>%
f^{\prime }%
&%
&%
\\ %
\newline
\end{CD}%
$ & .%
\end{tabular}%
\end{equation*}%
And furthermore, the following conditions are satisfied whenever the
compositions are defined%
\begin{gather*}
\xymatrix@R=1pc@C=0.1pc{
& 
f \otimes (g \otimes (h \otimes k))
\ar[rr]^{f \otimes \alpha_{g,h,k} }
\ar[ldd] _{ \alpha_{f,g,h\otimes k} }
& & 
f \otimes ((g \otimes h) \otimes k)
\ar[rdd]^{ \alpha_{f,g\otimes h,k} }
& \\
& & & & \\
(f \otimes g) \otimes (h \otimes k)   
\ar[rrdd]_{ \alpha_{f\otimes g,h,k} }
 &  &  &  &  
(f \otimes (g \otimes h)) \otimes k  
\ar[ddll]^{\alpha_{f,g,h} \otimes k}
\\
& & & & \\
 & & 
((f \otimes g) \otimes h) \otimes k
 &  &}
\\
\xymatrix@=4pc{
 f\otimes (1\otimes g)   \ar[rr]^{\alpha_{f,1,g}}
            \ar[rd]_{f\otimes \lambda_g}
&& (f\otimes 1)\otimes g   \ar[ld]^{\rho_f \otimes g} \\
& f\otimes g}%
\begin{array}[t]{c}
\\ 
.%
\end{array}%
\end{gather*}

Examples of pseudo-categories internal to Cat include the usual bicategories
of Spans, Bimodules, homotopies, ... where in each case it is also allowed
to consider the natural morphisms between the objects in order to obtain a
vertical categorical structure. For example in the case of spans\ we would
have sets as objects, maps as morphisms, spans $A\longleftarrow
S\longrightarrow B$ as pseudo-morphisms and the cells being triples $\left(
h,k,l\right) $ with the\ following two squares commutative%
\begin{equation*}
\begin{tabular}{ll}
$%
\begin{CD}%
\newline
A%
@<<<%
S%
@>%
>%
>%
B%
\\ %
\newline
@V%
h%
V%
V%
@V%
k%
V%
V%
@V%
V%
l%
V%
\\ %
\newline
A^{\prime }%
@<<<%
S^{\prime }%
@>%
>%
>%
B%
&%
&%
\\ %
\newline
\end{CD}%
$ & .%
\end{tabular}%
\end{equation*}

A pseudo-category in Cat has the following structures: a category (with
objects and morphisms); a category (with pseudo-morphisms and cells); a
bicategory (considering only the morphisms that are identities); a double
category (if all the special cells are identity cells).

Other examples as Cat (with modules as pseudo-morphisms) may be found in 
\cite{Leinster} or \cite{Pare}.

The present description of pseudo double category (internal pseudo-category
in Cat) is the same given by Leinster \cite{Leinster} and differs from the
one considered by Grandis and Par\'{e} \cite{Pare} in the sense that they
also have 
\begin{equation*}
id_{A}=id_{A}\otimes id_{A}.
\end{equation*}

In the following sections we will provide a complete description of
pseudo-functors, natural and pseudo-natural transformations and
pseudo-modifications. We prove that all the compositions are well defined
(except for the horizontal composition of pseudo-natural transformations
which is only defined up to an isomorphism). In the end we show that the
category of pseudo-categories (internal to some ambient \ 2-category \textbf{%
C}) is Cartesian closed up to isomorphism. We will give all the definitions
in terms of the internal structure to some ambient 2-category and also
explain what is obtained in the case where the ambient 2-category is Cat.
While doing some proofs we will \ make use of Yoneda embedding and consider
the diagrams in Cat rather than in the abstract ambient 2-category.

We will also freely use known definitions and results from \cite{ML},\cite{Benabou},\cite%
{Gray} and \cite{Power}.

\section{Pseudo-Functors}

Let \textbf{C} be a 2-category and suppose%
\begin{eqnarray}
C &=&\left( C_{0},C_{1},d,c,e,m,\alpha ,\lambda ,\rho \right) ,
\label{C and C'} \\
C^{\prime } &=&\left( C_{0}^{\prime },C_{1}^{\prime },d^{\prime },c^{\prime
},e^{\prime },m^{\prime },\alpha ^{\prime },\lambda ^{\prime },\rho ^{\prime
}\right)  \notag
\end{eqnarray}%
are two pseudo-categories in \textbf{C}.

A pseudo-functor $F:C\longrightarrow C^{\prime }$ is a system 
\begin{equation*}
F=\left( F_{0},F_{1},\mu ,\varepsilon \right)
\end{equation*}%
where $F_{0}:C_{0}\longrightarrow C_{0}^{\prime },F_{1}:C_{1}\longrightarrow
C_{1}^{\prime }$ are morphisms of \textbf{C}, 
\begin{equation*}
\mu :F_{1}m\longrightarrow m^{\prime }\left( F_{1}\times
_{F_{0}}F_{1}\right) \;,\;\varepsilon :F_{1}e\longrightarrow e^{\prime
}F_{0},
\end{equation*}%
are 2-cells of \textbf{C} (that are isomorphisms\footnote{%
Some authors (example Grandis and Par\'{e} in \cite{Pare,Pare2}) consider
the notion of pseudo - which corresponds to the present one - but also
consider the notions of lax and colax where the 2-cells may not be
isomorphisms.}), the following conditions are satisfied%
\begin{eqnarray}
d^{\prime }F_{1} &=&F_{0}d,  \label{dom cod for F0,F1} \\
c^{\prime }F_{1} &=&F_{0}c,  \notag
\end{eqnarray}%
\begin{eqnarray}
d^{\prime }\circ \mu &=&1_{F_{0}d\pi _{2}},  \label{dom cod for miuF} \\
c^{\prime }\circ \mu &=&1_{F_{0}c\pi _{1}},  \notag
\end{eqnarray}%
\begin{eqnarray}
d^{\prime }\circ \varepsilon &=&1_{F_{0}},  \label{dom cod for epsilonF} \\
c^{\prime }\circ \varepsilon &=&1_{F_{0}},  \notag
\end{eqnarray}%
and the following diagrams commute%
\begin{equation}
\xymatrix{
&
\bullet 
     \ar[r]^{\mu \left( 1\times _{C_{0}}m\right)} 
     \ar[ld]_{F_{1}\alpha}
&
\bullet
      \ar[rd]^{m^{\prime }\left( 1_{F_{1}}\times \mu \right) }
\\
\bullet
      \ar[rd]_{\mu \left( m\times_{C_{0}}1\right)}
 & & &
\bullet
      \ar[ld]^{ \alpha ^{\prime }\left( F_{1}\times _{F_{0}}F_{1}\times _{F_{0}}F_{1}\right) }
\\
&
\bullet
     \ar[r]_{m^{\prime }\left( \mu \times 1_{F_{1}}\right) }
& \bullet }%
\begin{array}[t]{c}
\\ 
,%
\end{array}
\label{hexagon F}
\end{equation}%
\begin{eqnarray}
&&%
\begin{tabular}{ll}
$%
\begin{CD}%
\newline
\bullet 
@>%
F_{1}\rho 
>%
>%
\bullet 
\\ %
\newline
@V%
\mu \left\langle 1,ed\right\rangle 
V%
V%
@A%
A%
\rho ^{\prime }F_{1}%
A%
\\ %
\newline
\bullet 
@>%
m^{\prime }\left( 1_{F_{1}}\times \varepsilon \right) 
>%
>%
\bullet 
&%
&%
\\ %
\newline
\end{CD}%
$ & \ \ \ \ \ \ \ \ ,%
\end{tabular}
\label{squares coherence for F} \\
&&%
\begin{tabular}{ll}
$%
\begin{CD}%
\newline
\bullet 
@>%
F_{1}\lambda 
>%
>%
\bullet 
\\ %
\newline
@V%
\mu \left\langle ec,1\right\rangle 
V%
V%
@A%
A%
\lambda ^{\prime }F_{1}%
A%
\\ %
\newline
\bullet 
@>%
m^{\prime }\left( \varepsilon \times 1_{F_{1}}\right) 
>%
>%
\bullet 
&%
&%
\\ %
\newline
\end{CD}%
$ & \ \ \ \ \ \ \ \ \ .%
\end{tabular}
\notag
\end{eqnarray}

Consider the particular case of \textbf{C}=Cat. Let 
\begin{equation*}
\xymatrix@R=1pc@C=5pc{
A \ar[r]|{f}="f" \ar[d]_{a} & B \ar[d]^{b} \\
A^{\prime} \ar[r]|{f^{\prime}}="k" & B^{\prime}  \\
\ar@{=>}^{\varphi} "f";"k"}%
\end{equation*}%
be a cell in the pseudo-category $C$. A pseudo-functor $F:C\longrightarrow
C^{\prime }$, consists of four maps (sending objects to objects, morphisms
to morphisms, pseudo-morphisms to pseudo-morphisms and cells to cells - that
we will denote only by $F$ to keep notation simple)%
\begin{equation*}
\xymatrix@R=1pc@C=5pc{
FA \ar[r]|{Ff}="f" \ar[d]_{Fa} & FB \ar[d]^{Fb} \\
FA^{\prime} \ar[r]|{Ff^{\prime}}="k" & FB^{\prime}  \\
\ar@{=>}^{F\varphi} "f";"k"}%
;
\end{equation*}%
a special cell $\mu _{f,g}$ 
\begin{equation*}
\squarecell{FA}   {F(g\otimes f)}    {FC}
                  {1}   {\mu_{f,g}}    {1}
                  {FA}   {Fg\otimes Ff}    {FC}%
\end{equation*}%
to each pair of composable pseudo-morphisms $f,g$; a special cell $%
\varepsilon _{A}$%
\begin{equation*}
\squarecell{FA}   {F(id_A)}    {FA}
                  {1}   {\epsilon_A}    {1}
                  {FA}   {id_{FA}}    {FA}%
\end{equation*}%
to each object $A$, and satisfying the commutativity of the following
diagrams%
\begin{equation*}
\hspace{-1cm}\hspace{-1cm}%
\xymatrix{
&
F\left( f\otimes \left( g\otimes h\right) \right) 
      \ar[r]^{F\left( \alpha _{f,g,h}\right) }
      \ar[ld]_{\mu _{f,g\otimes h}}
&
F\left( \left( f\otimes g\right) \otimes h\right) 
       \ar[rd]^{\mu _{f\otimes g,h}}
\\
F\left( f\right) \otimes F\left( g\otimes h\right) 
        \ar[rd]_{F\left( f\right) \otimes \mu _{g,h}}
& & &
F\left( f\otimes g\right) \otimes F\left( h\right) 
          \ar[ld]^{\mu _{f,g}\otimes F\left( h\right) }
\\
&
F\left( f\right) \otimes \left( F\left( g\right) \otimes F\left( h\right) \right) 
           \ar[r]_{\alpha _{Ff,Fg,Fh}^{\prime }}
&
\left( F\left( f\right) \otimes F\left( g\right) \right) \otimes F\left( h\right)
}%
\begin{array}[t]{c}
\\ 
\\ 
,%
\end{array}%
\end{equation*}%
\begin{equation*}
\begin{tabular}{ll}
$%
\begin{CD}%
\newline
F\left( f\otimes id_{A}\right) 
@>%
F\left( \rho _{f}\right) 
>%
>%
F\left( f\right) 
\\ %
\newline
@V%
\mu _{f,id_{A}}%
V%
V%
@A%
A%
\rho _{Ff}^{\prime }%
A%
\\ %
\newline
F\left( f\right) \otimes F\left( id_{A}\right) 
@>%
F\left( f\right) \otimes \varepsilon _{A}%
>%
>%
F\left( f\right) \otimes id_{F\left( A\right) }%
&%
&%
\\ %
\newline
\end{CD}%
$ & \ \ \ \ \ \ \ \ ,%
\end{tabular}%
\end{equation*}%
\begin{equation*}
\begin{tabular}{ll}
$%
\begin{CD}%
\newline
F\left( id_{B}\otimes f\right) 
@>%
F\left( \lambda _{f}\right) 
>%
>%
F\left( f\right) 
\\ %
\newline
@V%
\mu _{id_{B},f}%
V%
V%
@A%
A%
\lambda _{Ff}^{\prime }%
A%
\\ %
\newline
F\left( id_{B}\right) \otimes F\left( f\right) 
@>%
\varepsilon _{B}\otimes F\left( f\right) 
>%
>%
id_{F\left( B\right) }\otimes F\left( f\right) 
&%
&%
\\ %
\newline
\end{CD}%
$ & \ \ \ \ \ \ \ \ ,%
\end{tabular}%
\end{equation*}%
whenever the pseudo-compositions are defined.

\bigskip

Return to the general case.

Let $F:C\longrightarrow C^{\prime }$ and $G:C^{\prime }\longrightarrow
C^{\prime \prime }$ be pseudo-functors in a 2-category \textbf{C}. Consider $%
C$ and $C^{\prime }$ as in $\left( \ref{C and C'}\right) $ and let%
\begin{eqnarray*}
C^{\prime \prime } &=&\left( C_{0}^{\prime \prime },C_{1}^{\prime \prime
},d^{\prime \prime },c^{\prime \prime },e^{\prime \prime },m^{\prime \prime
},\alpha ^{\prime \prime },\lambda ^{\prime \prime },\rho ^{\prime \prime
}\right) , \\
F &=&\left( F_{0},F_{1},\mu ^{F},\varepsilon ^{F}\right) , \\
G &=&\left( G_{0},G_{1},\mu ^{G},\varepsilon ^{G}\right) .
\end{eqnarray*}%
The composition of the pseudo-functors $F$ and $G$ is defined by the formula%
\begin{equation}
GF=\left( G_{0}F_{0},G_{1}F_{1},\left( \mu ^{G}\circ \left( F_{1}\times
_{F_{0}}F_{1}\right) \right) \cdot \left( G_{1}\circ \mu ^{F}\right) ,\left(
\varepsilon ^{G}\circ F_{0}\right) \cdot \left( G_{1}\circ \varepsilon
^{F}\right) \right)  \label{GF}
\end{equation}%
where $\circ $ represents the horizontal composition in \textbf{C} and $%
\cdot $ represents the vertical composition, as displayed in the diagram
below%
\begin{equation*}
\begin{tabular}{ll}
$%
\begin{CD}%
\newline
C_{1}\times _{C_{0}}C_{1}%
@>%
m%
>%
>%
C_{1}%
@<%
e%
<%
<%
C_{0}%
\\ %
\newline
@V%
F_{1}\times _{F_{0}}F_{1}%
V%
V%
\hspace{-1.5cm}\mu ^{F}\Downarrow \ \ \ \ \ 
@V%
F_{1}%
V%
V%
\hspace{-1cm}\varepsilon ^{F}\Downarrow \ \ \ 
@V%
V%
F_{0}%
V%
\\ %
\newline
C_{1}^{\prime }\times _{C_{0}^{\prime }}C_{1}^{\prime }%
@>%
m^{\prime }%
>%
>%
C_{1}^{\prime }%
@<%
e^{\prime }%
<%
<%
C_{0}^{\prime }%
\\ %
$\newline
$%
@V%
G_{1}\times _{G_{0}}G_{1}%
V%
V%
\hspace{-1.5cm}\mu ^{G}\Downarrow \ \ \ \ \ 
@V%
G_{1}%
V%
V%
\hspace{-1cm}\varepsilon ^{G}\Downarrow \ \ \ 
@V%
V%
G_{0}%
V%
\\ %
\newline
C_{1}^{\prime \prime }\times _{C_{0}^{\prime \prime }}C_{1}^{\prime \prime }%
@>%
m^{\prime \prime }%
>%
>%
C_{1}^{\prime \prime }%
@<%
e^{\prime \prime }%
<%
<%
C_{0}^{\prime \prime }%
\\ %
\newline
\end{CD}%
$ & .%
\end{tabular}%
\end{equation*}

\begin{proposition}
The above formula to compose pseudo-functors is well defined.
\end{proposition}

\begin{proof}
Consider the system%
\begin{equation*}
GF=\left( G_{0}F_{0},G_{1}F_{1},\mu ^{GF},\varepsilon ^{GF}\right)
\end{equation*}%
with $\mu ^{GF},\varepsilon ^{GF}$ as in $\left( \ref{GF}\right) $. We will
show that $GF$ is a pseudo-functor from the pseudo-category $C$ to the
pseudo-category $C^{\prime \prime }$.

It is clear that $G_{0}F_{0}:C_{0}\longrightarrow C_{0}^{\prime \prime
},G_{1}F_{1}:C_{1}\longrightarrow C_{1}^{\prime \prime },$ are morphisms of
the ambient 2-category \textbf{C} and$\ \mu ^{GF}:G_{1}F_{1}m\longrightarrow
m^{\prime \prime }\left( G_{1}F_{1}\times _{G_{0}F_{0}}G_{1}F_{1}\right) ,\
\varepsilon ^{GF}:G_{1}F_{1}e\longrightarrow e^{\prime \prime }G_{0}F_{0}$
are 2-cells of \textbf{C}\ and they are isomorphisms.

Conditions $\left( \ref{dom cod for F0,F1}\right) $ are satisfied and 
\begin{eqnarray*}
d^{\prime \prime }\mu ^{GF} &=&d\left( \left( \mu ^{G}\circ \left(
F_{1}\times _{F_{0}}F_{1}\right) \right) \cdot \left( G_{1}\circ \mu
^{F}\right) \right) \\
&=&\left( d\circ \mu ^{G}\circ \left( F_{1}\times _{F_{0}}F_{1}\right)
\right) \cdot \left( d\circ G_{1}\circ \mu ^{F}\right) \\
&=&\left( 1_{G_{0}d^{\prime }\pi _{2}^{\prime }}\circ \left( F_{1}\times
_{F_{0}}F_{1}\right) \right) \cdot \left( G_{0}\circ d^{\prime }\circ \mu
^{F}\right) \\
&=&\left( 1_{G_{0}d^{\prime }\pi _{2}^{\prime }\left( F_{1}\times
_{F_{0}}F_{1}\right) }\right) \cdot \left( G_{0}\circ 1_{F_{0}d\pi
_{2}}\right) \\
&=&1_{G_{0}d^{\prime }F_{1}\pi _{2}}\cdot 1_{G_{0}F_{0}d\pi _{2}} \\
&=&1_{G_{0}F_{0}d\pi _{2}},
\end{eqnarray*}%
as well $c^{\prime \prime }\mu ^{GF}=1_{G_{0}F_{0}c\pi _{1}},$ hence $\left( %
\ref{dom cod for miuF}\right) $ holds. Also%
\begin{eqnarray*}
d^{\prime \prime }\varepsilon ^{GF} &=&d^{\prime \prime }\left( \left(
\varepsilon ^{G}\circ F_{0}\right) \cdot \left( G_{1}\circ \varepsilon
^{F}\right) \right) \\
&=&\left( d^{\prime \prime }\circ \varepsilon ^{G}\circ F_{0}\right) \cdot
\left( d^{\prime \prime }G_{1}\circ \varepsilon ^{F}\right) \\
&=&\left( 1_{G_{0}}\circ F_{0}\right) \cdot \left( G_{0}d^{\prime }\circ
\varepsilon ^{F}\right) \\
&=&1_{G_{0}F_{0}}\cdot \left( G_{0}\circ 1_{F_{0}}\right) \\
&=&1_{G_{0}F_{0}}\cdot 1_{G_{0}F_{0}}=1_{G_{0}F_{0}},
\end{eqnarray*}%
and similarly $c^{\prime \prime }\varepsilon ^{GF}=1_{G_{0}F_{0}}$, so
conditions $\left( \ref{dom cod for epsilonF}\right) $ are satisfied.

Commutativity of diagrams $\left( \ref{hexagon F}\right) ,$ $\left( \ref%
{squares coherence for F}\right) $ follows from Yoneda Lemma and the
commutativity of the following diagrams%
\begin{eqnarray*}
&&\hspace{-1cm}\hspace{-1cm}\hspace{-1cm}%
\xymatrix{
&   
    H_{f \otimes (g\otimes h)}             \ar[rrr]^{H(\alpha_{f,g,h})}
          \ar[d]^{}
          \ar[ldd]_{\mu^{H}_{f , g\otimes h}}
& & &
    H_{(f\otimes g)\otimes h)}            \ar[rdd]^{\mu^{H}_{f\otimes g,h}}
          \ar[d]^{}   
\\
&
    G_{F_f\otimes F_{g\otimes h}}            \ar[ld]^{}
          \ar[rd]^{}
& & &
    G_{F_{f\otimes g} \otimes F_h}             \ar[ld]^{}
          \ar[rd]^{}
\\
    H_f\otimes H_{g\otimes h}            \ar[rd]^{}
          \ar[rdd]_{H_f \otimes \mu^{H}_{g,h}}
& &
     (1)             \ar[ld]^{}
          \ar[r]^{G(\alpha^{\prime}_{Ff,Fg,Fh})}
&
      (2)             \ar[rd]^{}
& & 
    H_{f\otimes g} \otimes H_h            \ar[ld]^{}
          \ar[ldd]^{\mu^{H}_{f,g}\otimes H_h}
\\
& 
    F_f \otimes G_{F_g\otimes F_h}             \ar[d]^{}
& & &
     G_{F_f\otimes F_g} \otimes H_h             \ar[d]^{}
\\
& 
    H_f \otimes (H_g\otimes H_h)                 \ar[rrr]^{\alpha^{\prime\prime}_{Hf,Hg,Hh}}
& & &     ( H_f \otimes H_g) \otimes H_h      \\}%
\begin{array}[t]{c}
\\ 
\\ 
\\ 
\\ 
\end{array}
\\
&&%
\xymatrix@=1pc{
  H_{id_B \otimes f}
            \ar[rrr]^{H(\lambda_f)}
	\ar[rd]
	\ar[ddd]_{\mu^{H}_{id_B,f}}
& & &
  H_f
\\
& 
  G_{F_{id_B} \otimes F_f}
	\ar[r]
	\ar[ldd]
&
  G_{id_{F_B} \otimes F_f}
	\ar[d]
\\
& &
  G_{id_{F_B}} \otimes H_f
	\ar[dr]
\\
  H_{id_B} \otimes H_f
	\ar[urr]
	\ar[rrr]_{\epsilon_B^{H}\otimes Hf}
& & &
  id_{H_B} \otimes H_f
	\ar[uuu]_{\lambda^{\prime\prime}_{Hf}}
\\}%
\begin{array}[t]{c}
\\ 
\\ 
\end{array}
\\
&&%
\xymatrix@=1pc{
  H_{f \otimes id_A}
            \ar[rrr]^{H(\rho_f)}
	\ar[rd]
	\ar[ddd]_{\mu^{H}_{f,id_A}}
& & &
  H_f
\\
& 
  G_{F_{f} \otimes F_{id_A}}
	\ar[r]
	\ar[ldd]
&
  G_{F_f \otimes id_{F_A}}
	\ar[d]
\\
& &
  H_f \otimes G_{id_{F_A}}
	\ar[dr]
\\
  H_{f} \otimes H_{id_A}
	\ar[urr]
	\ar[rrr]_{Hf\otimes \epsilon_A^{H}}
& & &
  F_f \otimes id_{H_A}
	\ar[uuu]_{\rho^{\prime\prime}_{Hf}}
\\}%
\begin{array}[t]{c}
\\ 
\\ 
\end{array}%
\end{eqnarray*}%
where $(1)=G_{F_{f}\otimes (F_{g}\otimes F_{h})}$ and $(2)=G_{(F_{f}\otimes
F_{g})\otimes F_{h}}$. We also use the abbreviations $H=GF$ and $F_{f}$ or $%
Ff$ instead of $F\left( f\right) $ to save space in the diagram.
\end{proof}

Composition of pseudo-functors is associative and there is an identity
pseudo-functor for every pseudo-category, namely the pseudo-functor%
\begin{equation*}
1_{C}=\left( 1_{C_{0}},1_{C_{1}},1_{m},1_{e}\right)
\end{equation*}%
for the pseudo-category%
\begin{equation*}
C=\left( C_{0},C_{1},d,c,e,m,\alpha ,\lambda ,\rho \right) .
\end{equation*}

Given a 2-category \textbf{C}, we define the category PsCat(\textbf{C})
consisting of all pseudo-categories and pseudo-functors internal to \textbf{C%
}.

\section{Natural and pseudo-natural transformations}

Let \textbf{C} be a 2-category and suppose 
\begin{eqnarray}
C &=&\left( C_{0},C_{1},d,c,e,m,\alpha ,\lambda ,\rho \right) ,
\label{C and C' again} \\
C^{\prime } &=&\left( C_{0}^{\prime },C_{1}^{\prime },d^{\prime },c^{\prime
},e^{\prime },m^{\prime },\alpha ^{\prime },\lambda ^{\prime },\rho ^{\prime
}\right)  \notag
\end{eqnarray}%
are pseudo-categories in \textbf{C} and 
\begin{eqnarray}
F &=&\left( F_{0},F_{1},\mu ^{F},\varepsilon ^{F}\right) ,  \label{F and G}
\\
G &=&\left( G_{0},G_{1},\mu ^{G},\varepsilon ^{G}\right)  \notag
\end{eqnarray}%
are pseudo-functors from $C$ to $C^{\prime }$.

A \textbf{natural transformation} $\theta :F\longrightarrow G$ is a pair $%
\theta =\left( \theta _{0},\theta _{1}\right) $ of 2-cells of \textbf{C}%
\begin{eqnarray*}
\theta _{0} &:&F_{0}\longrightarrow G_{0} \\
\theta _{1} &:&F_{1}\longrightarrow G_{1}
\end{eqnarray*}%
satisfying%
\begin{eqnarray*}
d^{\prime }\circ \theta _{1} &=&\theta _{0}\circ d \\
c^{\prime }\circ \theta _{1} &=&\theta _{0}\circ c
\end{eqnarray*}%
and the commutativity of the following diagrams of 2-cells%
\begin{eqnarray*}
&&%
\begin{tabular}{ll}
$%
\begin{CD}%
\newline
\bullet 
@>%
\theta _{1}\circ m%
>%
>%
\bullet 
\\ %
\newline
@V%
\mu ^{F}%
V%
V%
@V%
V%
\mu ^{G}%
V%
\\ %
\newline
\bullet 
@>%
m^{\prime }\circ \left( \theta _{1}\times _{\theta _{0}}\theta _{1}\right) 
>%
>%
\bullet 
&%
&%
\\ %
\newline
\end{CD}%
$ & 
\end{tabular}
\\
&&%
\begin{tabular}{ll}
$%
\begin{CD}%
\newline
\bullet 
@>%
\theta _{1}\circ e%
>%
>%
\bullet 
\\ %
\newline
@V%
\varepsilon ^{F}%
V%
V%
@V%
V%
\varepsilon ^{G}%
V%
\\ %
\newline
\bullet 
@>%
e^{\prime }\circ \theta _{0}%
>%
>%
\bullet 
&%
&%
\\ %
\newline
\end{CD}%
$ & .%
\end{tabular}%
\end{eqnarray*}

A \textbf{pseudo-natural transformation} $T:F\longrightarrow G$ is a pair 
\begin{equation*}
T=\left( t,\tau \right)
\end{equation*}%
where $t:C_{0}\longrightarrow C_{1}^{\prime }$ is a morphism of \textbf{C}, 
\begin{equation*}
\tau :m^{\prime }\left\langle G_{1},td\right\rangle \longrightarrow
m^{\prime }\left\langle tc,F_{1}\right\rangle
\end{equation*}%
is a 2-cell (that is an isomorphism); the following conditions are satisfied%
\begin{eqnarray}
d^{\prime }t &=&F_{0}  \label{dom cod t} \\
c^{\prime }t &=&G_{0}  \notag
\end{eqnarray}%
\begin{eqnarray}
d^{\prime }\circ \tau &=&1_{d^{\prime }F_{1}}  \label{dom cod tau} \\
c^{\prime }\circ \tau &=&1_{c^{\prime }G_{1}}  \notag
\end{eqnarray}%
and the following diagrams of 2-cells are commutative$\footnote{$G_{1}\times
_{G_{0}}t\times _{F_{0}}F_{1}:C_{1}\times _{C_{0}}C_{0}\times
_{C_{0}}C_{1}\longrightarrow C_{1}\times _{C_{0}}C_{1}\times _{C_{0}}C_{1}$%
\par
$t\times _{F_{0}}F_{1}\times _{F_{0}}F_{1}:C_{0}\times _{C_{0}}C_{1}\times
_{C_{0}}C_{1}\longrightarrow C_{1}\times _{C_{0}}C_{1}\times _{C_{0}}C_{1}$%
\par
$G_{1}\times _{G_{0}}G_{1}\times _{G_{0}}t:C_{1}\times _{C_{0}}C_{1}\times
_{C_{0}}C_{0}\longrightarrow C_{1}\times _{C_{0}}C_{1}\times _{C_{0}}C_{1}$}$%
\begin{equation}
\xymatrix{
&
\bullet 
      \ar[r]^{\alpha ^{-1}\left( G_{1}\times _{G_{0}}G_{1}\times _{G_{0}}t\right)}
      \ar[ld]_{m^{\prime }\left\langle \mu _{G}^{-1},td\pi _{2}\right\rangle }
&
 \bullet
        \ar[rd]^{ m^{\prime }\left\langle G_{1}\pi _{1},\tau \pi _{2}\right\rangle }
\\
\bullet 
          \ar[d]_{ \tau m}
& & &
\bullet
           \ar[d]^{\alpha \left( G_{1}\times _{G_{0}}t\times _{F_{0}}F_{1}\right) }
\\
\bullet
            \ar[rd]_{ m^{\prime }\left\langle {tc\pi _{1}},\mu _{F}\right\rangle}
& & &
 \bullet
             \ar[ld]^{m^{\prime }\left\langle \tau \pi _{1},F_{1}\pi _{2}\right\rangle }
\\
&
\bullet
             \ar[r]_{\alpha \left( t\times _{F_{0}}F_{1}\times _{F_{0}}F_{1}\right) }
&
\bullet }%
\begin{array}[t]{c}
\\ 
,%
\end{array}
\label{Octogon}
\end{equation}%
\begin{equation}
\xymatrix{
&
\bullet 
          \ar[rr]^{\tau e}
          \ar[dl]_{\lambda ^{\prime }t}
& &
\bullet 
           \ar[dr]^{m^{\prime }\left\langle 1_{t},\varepsilon _{F}\right\rangle }
\\
\bullet
            \ar[rrd]_{m^{\prime }\left\langle \varepsilon _{G},1_{t}\right\rangle}
& & & & 
 \bullet 
          \ar[lld]^{\rho ^{\prime }t}
\\
& &
\bullet }%
\begin{array}[t]{c}
\\ 
.%
\end{array}
\label{pentagon tau}
\end{equation}

\medskip In the case \textbf{C}=Cat: let $W,W^{\prime }$ be two
pseudo-categories in Cat, and $F,G:W\longrightarrow W^{\prime }$ two
pseudo-functors. Given a cell%
\begin{equation*}
\squarecell{A}   {f}    {B}
                  {a}   {\varphi}    {b}
                  {A^\prime}   {f^\prime}    {B^\prime}%
\end{equation*}%
in $W$, we will write%
\begin{equation*}
\squarecell{FA}   {Ff}    {FB}
                  {Fa}   {F\varphi}    {Fb}
                  {FA^\prime}   {Ff^\prime}    {FB^\prime}%
\begin{array}[t]{c}
\\ 
\ \ \text{and\ \ }%
\end{array}%
\squarecell{GA}   {Gf}    {GB}
                  {Ga}   {G\varphi}    {Gb}
                  {GA^\prime}   {Gf^\prime}    {GB^\prime}%
\end{equation*}%
for the image of $\varphi $ under $F$ and $G$.

The description of natural and pseudo-natural transformations in this
particular case is as follows:

- While a natural transformation $\theta :F\longrightarrow G$ is a family of
cells%
\begin{equation*}
\squarecell{FA}   {Ff}    {FB}
                  {\theta_A}   {\theta_f}    {\theta_B}
                  {GA}   {Gf}    {GB}%
\begin{array}[t]{c}
\\ 
,%
\end{array}%
\end{equation*}%
one for each pseudo-morphism $f$ in $W$, such that for every cell $\varphi $
in $W$, the square%
\begin{equation*}
\begin{tabular}{ll}
$%
\begin{CD}%
\newline
Ff%
@>%
\theta _{f}%
>%
>%
Gf%
\\ %
\newline
@V%
F\varphi 
V%
V%
@V%
V%
G\varphi 
V%
\\ %
\newline
Ff^{\prime }%
@>%
\theta _{f^{\prime }}%
>%
>%
Gf^{\prime }%
&%
&%
\\ %
\newline
\end{CD}%
$ & 
\end{tabular}%
\end{equation*}%
is commutative as displayed in the picture below%
\begin{equation*}
\xymatrix@R=1pc@C=5pc{
	\ar[rr]|{Ff}="Ff"
	\ar[d]_{Fa}
	\ar[rdd]_{\theta_A}
& &
	\ar[d]^{Fb}
	\ar[rdd]^{\theta_B}
\\
	\ar[rr]|{Ff^{\prime}}="Ffp"
	\ar[rdd]_{\theta_{A^\prime}}
& &
	\ar[rdd]^{\theta_{B^\prime}}
\\
&
	\ar[rr]|{Gf}="Gf"
	\ar[d]_{Ga}
& &
	\ar[d]^{Gb}
\\
&
	\ar[rr]|{Gf^\prime}="Gfp"
& & 
\\
\ar@{=>}^{\theta_f} "Ff";"Gf" 
\ar@{=>}^{\theta_{f^\prime}} "Ffp";"Gfp" 
\ar@{=>}_{F\varphi} "Ff";"Ffp" 
\ar@{=>}^{G\varphi} "Gf";"Gfp" 
}%
\begin{array}[t]{c}
\\ 
;%
\end{array}%
\end{equation*}%
and furthermore, given two composable pseudo-morphisms $g,f$ and an object $%
A $ in $W$, the following squares are commutative%
\begin{equation*}
\begin{tabular}{ll}
$%
\begin{CD}%
\newline
F\left( g\otimes f\right) 
@>%
\mu _{g,f}^{F}%
>%
>%
Fg\otimes Ff%
\\ %
\newline
@V%
\theta _{g\otimes f}%
V%
V%
@V%
V%
\theta _{g}\otimes \theta _{f}%
V%
\\ %
\newline
G\left( g\otimes f\right) 
@>%
\mu _{g,f}^{G}%
>%
>%
Gg\otimes Gf%
&%
&%
\\ %
\newline
\end{CD}%
$ & 
\end{tabular}%
\end{equation*}%
\begin{equation*}
\begin{tabular}{ll}
$%
\begin{CD}%
\newline
F\left( id_{A}\right) 
@>%
\varepsilon _{A}^{F}%
>%
>%
id_{FA}%
\\ %
\newline
@V%
\theta _{id_{A}}%
V%
V%
@V%
V%
id_{\theta _{A}}%
V%
\\ %
\newline
G\left( id_{A}\right) 
@>%
>%
\varepsilon _{A}^{G}%
>%
id_{GA}%
&%
&%
\\ %
\newline
\end{CD}%
$ & .%
\end{tabular}%
\end{equation*}%
\medskip - Rather a pseudo-natural transformation $T:F\longrightarrow G$
consists of two families of cells%
\begin{equation*}
\squarecell{FA}   {t_A}    {GA}
                  {Fa}   {t_a}    {Ga}
                  {FA^\prime}   {t_{A^\prime}}    {GA^\prime}%
\end{equation*}%
and%
\begin{equation*}
\squarecell{FA}   {Gf\otimes t_A}    {GB}
                  {1}   {\tau_f}    {1}
                  {FA}   {t_B\otimes Ff}    {GB}%
\end{equation*}%
with $a$ a morphism and $f$ a pseudo-morphism of $W$, as displayed in the
following picture%
\begin{equation*}
\xymatrix@R=1pc@C=5pc{
	\ar[rr]|{Ff}="Ff"
	\ar[d]_{Fa}
	\ar[rdd]|{t_A}="t_A"
& &
	\ar[d]^{Fb}
	\ar[rdd]|{t_B}="t_B"
\\
	\ar[rr]|{Ff^{\prime}}="Ffp"
	\ar[rdd]|{t_{A^\prime}}="t_Ap"
& &
	\ar[rdd]|{t_{B^\prime}}="t_Bp"
\\
&
	\ar[rr]|{Gf}="Gf"
	\ar[d]_{Ga}
	\ar@{=>}[uur]_{\tau_f}
& &
	\ar[d]^{Gb}
\\
&
	\ar[rr]|{Gf^\prime}="Gfp"
	\ar@{=>}[uur]_{\tau_{f^\prime}}
& & 
\\
\ar@{=>}^{t_a} "t_A";"t_Ap" 
\ar@{=>}^{t_b} "t_B";"t_Bp" 
\ar@{=>}_{F\varphi} "Ff";"Ffp" 
\ar@{=>}^{G\varphi} "Gf";"Gfp" 
}%
\end{equation*}%
such that ($t$ is a functor)%
\begin{eqnarray*}
t_{a^{\prime }a} &=&t_{a}t_{a^{\prime }} \\
t_{1_{A}} &=&1_{t_{A}}
\end{eqnarray*}%
($\tau $ is natural)%
\begin{equation*}
\begin{tabular}{ll}
$%
\begin{CD}%
\newline
Gf\otimes t_{A}%
@>%
\tau _{f}%
>%
>%
t_{B}\otimes Ff%
\\ %
\newline
@V%
G\varphi \otimes t_{a}%
V%
V%
@V%
V%
t_{b}\otimes F\varphi 
V%
\\ %
\newline
Gf^{\prime }\otimes t_{A^{\prime }}%
@>%
\tau _{f^{\prime }}%
>%
>%
t_{B^{\prime }}\otimes Ff^{\prime }%
&%
&%
\\ %
\newline
\end{CD}%
$ & ,%
\end{tabular}%
\end{equation*}%
and for every two composable pseudo-morphisms $%
\xymatrix{A \ar[r]|{f} & B \ar[r]|{g} & C}%
$ , the following diagrams of cells in $W^{\prime }$ are commutative%
\begin{equation*}
\xymatrix{
&
\left( G\left( g\right) \otimes G\left( f\right) \right) \otimes t_{A}
               \ar[r]
               \ar[ld]
&
G\left( g\right) \otimes \left( G\left( f\right) \otimes t_{A}\right) 
               \ar[rd]
&\\
G\left( g\otimes f\right) \otimes t_{A}
               \ar[d]
& & &
G\left( g\right) \otimes \left( t_{B}\otimes F\left( f\right) \right) 
               \ar[d]
\\
t_{C}\otimes F\left( g\otimes f\right) 
               \ar[rd]
& & &
\left( G\left( g\right) \otimes t_{B}\right) \otimes F\left( f\right) 
               \ar[ld]
&\\
&
t_{C}\otimes \left( F\left( g\right) \otimes F\left( f\right) \right) 
               \ar[r]
&
\left( t_{C}\otimes F\left( g\right) \right) \otimes F\left( f\right) 
}%
\end{equation*}%
\begin{equation*}
\xymatrix{
&
G\left( id_{A}\right) \otimes t_{A}
             \ar[rr]
             \ar[ld]
& &
t_{A}\otimes F\left( id_{A}\right) 
             \ar[rd]
&\\
id_{GA}\otimes t_{A}
             \ar[rrd]
& & & &
t_{A}\otimes id_{FA}
             \ar[lld]
\\
& &
t_{A}
}%
\begin{array}[t]{c}
\\ 
.%
\end{array}%
\end{equation*}

\bigskip

Return to the general case.

Let \textbf{C} be a 2-category and suppose $C,C^{\prime },C^{\prime \prime }$
are pseudo-categories in \textbf{C}\ and $F,G,H:C\longrightarrow C^{\prime
}, $ $F^{\prime },G^{\prime }:C^{\prime }\longrightarrow C^{\prime \prime }$
are pseudo-functors. Natural transformations $\theta ,\theta ^{\prime },\dot{%
\theta}$%
\begin{equation*}
\xymatrix{
C
	\ar@/^1pc/[r]^{F}="F"
            \ar[r]^{\hspace{-0.5cm}G}
	\ar@{}[r]^{\,}="Gup"               \ar@{}[r]_{\,}="Gdown"   	\ar@/_1pc/[r]_{H}="H"
& C^\prime
	\ar@/^/[r]^{F^\prime}="Fp"
	\ar@/_/[r]_{G^\prime}="Gp"
& C^{\prime\prime}
\\
\ar@{=>}^{\theta} "F";"Gdown"
\ar@{=>}^{\dot{\theta}} "Gup";"H"
\ar@{=>}^{\theta^\prime} "Fp";"Gp"
}%
\end{equation*}%
may be composed horizontally with $\theta ^{\prime }\circ \theta =\left(
\theta _{0}^{\prime },\theta _{1}^{\prime }\right) \circ \left( \theta
_{0},\theta _{1}\right) =\left( \theta _{0}^{\prime }\circ \theta
_{0},\theta _{1}^{\prime }\circ \theta _{1}\right) $ obtained from the
horizontal composition of 2-cells of \textbf{C}, and vertically with $\dot{%
\theta}\cdot \theta =\left( \dot{\theta}_{0},\dot{\theta}_{1}\right) \cdot
\left( \theta _{0},\theta _{1}\right) =\left( \dot{\theta}_{0}\cdot \theta
_{0},\dot{\theta}_{1}\cdot \theta _{1}\right) $ obtained from the vertical
composition of 2-cells of \textbf{C}. Clearly both compositions are well
defined, are associative, have identities and satisfy the middle interchange
law. This fact may be stated as in the following theorem.

\begin{theorem}
\label{Theorem WCat is a 2Cat}Let \textbf{C} be a 2-category. The category
PsCat(\textbf{C}) (with pseudo-categories, pseudo-functors and natural
transformations) is a 2-category.
\end{theorem}

Composition of pseudo-natural transformations is much more delicate.

Again let \textbf{C} be a 2-category and suppose \ $C,C^{\prime }$ are
pseudo-categories in \textbf{C}, $F,G,H:C\longrightarrow C^{\prime }$ are
pseudo-functors (as above) and consider the pseudo-natural transformations%
\begin{equation*}
F\overset{T}{\longrightarrow }G\overset{S}{\longrightarrow }H
\end{equation*}%
with 
\begin{equation*}
T=\left( t,\tau \right) ,\ \ \ S=\left( s,\sigma \right) .
\end{equation*}%
\textbf{Vertical composition} of pseudo-natural transformations $S$ and $T$
is defined as%
\begin{equation}
S\otimes T=\left( m^{\prime }\left\langle s,t\right\rangle ,\sigma \otimes
\tau \right)  \label{ST}
\end{equation}%
where%
\begin{equation}
\sigma \otimes \tau =\alpha \left\langle sc,tc,F_{1}\right\rangle \cdot
m^{\prime }\left\langle 1_{sc},\tau \right\rangle \cdot \alpha
^{-1}\left\langle sc,G_{1},td\right\rangle \cdot m^{\prime }\left\langle
\sigma ,1_{td}\right\rangle \cdot \alpha \left\langle
H_{1},sd,td\right\rangle .  \label{tau_sigma}
\end{equation}

The above formula in the case \textbf{C}=Cat is expressed as follows%
\begin{equation*}
\begin{array}[t]{c}
\\ 
\left( s\otimes t\right) _{a}=s_{a}\otimes t_{a},\ \ \ \ \ \ \ 
\end{array}%
\xymatrix@R=1pc@C=5pc{
FA \ar[r]|{t_A}="f" \ar[d]_{Fa} & GA \ar[d]^{Ga} \ar[r]|{s_A}="g" & HA \ar[d]^{Ha} \\
FA^{\prime} \ar[r]|{t_{A^\prime}}="fp" & GA^{\prime}  \ar[r]|{s_{A^\prime}}="gp" & HA^{\prime}\\
\ar@{=>}^{t_a} "f";"fp"
\ar@{=>}^{s_a} "g";"gp"}%
\begin{array}[t]{c}
\\ 
;%
\end{array}%
\end{equation*}%
and%
\begin{equation*}
\left( \sigma \otimes \tau \right) _{f}=\alpha \left( s_{B}\otimes \tau
_{f}\right) \alpha ^{-1}\left( \sigma _{f}\otimes t_{A}\right) \alpha ,
\end{equation*}%
as displayed in the following picture%
\begin{equation*}
\begin{tabular}{ll}
$%
\begin{CD}%
\newline
Hf\otimes \left( s_{A}\otimes t_{A}\right) 
@>%
\left( \sigma \tau \right) _{f}%
>%
>%
\left( s_{B}\otimes t_{B}\right) \otimes Ff%
\\ %
\newline
@V%
\alpha 
V%
V%
@A%
A%
\alpha 
A%
\\ %
\newline
\left( Hf\otimes s_{A}\right) \otimes t_{A}%
&%
&%
s_{B}\otimes \left( t_{B}\otimes Ff\right) 
\\ %
\newline
@V%
\sigma _{f}\otimes t_{A}%
V%
V%
@A%
A%
s_{B}\otimes \tau _{f}%
A%
\\ %
\newline
\left( s_{B}\otimes Gf\right) \otimes t_{A}%
@>%
>%
\alpha ^{-1}%
>%
s_{B}\otimes \left( Gf\otimes t_{A}\right) 
&%
&%
\\ %
\newline
\end{CD}%
$ & \ \ \ \ \ \ \ \ .%
\end{tabular}%
\end{equation*}

Return to the general case.

\begin{theorem}
The vertical composition of pseudo-natural transformations is well defined.
\end{theorem}

\begin{proof}
Consider $C,C^{\prime }$ as in $\left( \ref{C and C' again}\right) $, $F,G$
as in $\left( \ref{F and G}\right) ,$ $H=\left( H_{0},H_{1},\mu
_{H},\varepsilon _{H}\right) $ and $S,T$ as above. Clearly $\left( st\right)
=m^{\prime }\left\langle s,t\right\rangle :C_{0}\longrightarrow
C_{1}^{\prime }$ is a morphism of \textbf{C} and $\sigma \tau :m^{\prime
}\left\langle H_{1},\left( st\right) d\right\rangle \longrightarrow
m^{\prime }\left\langle \left( st\right) c,F_{1}\right\rangle $ is a 2-cell
of \textbf{C} that is an isomorphism (is defined as a composition of
isomorphisms).

Conditions $\left( \ref{dom cod t}\right) $ and $\left( \ref{dom cod tau}%
\right) $ are satisfied%
\begin{eqnarray*}
d^{\prime }m^{\prime }\left\langle s,t\right\rangle &=&d^{\prime }\pi
_{2}^{\prime }\left\langle s,t\right\rangle \\
&=&d^{\prime }t \\
&=&F_{0},
\end{eqnarray*}%
also $c^{\prime }m^{\prime }\left\langle s,t\right\rangle =c^{\prime
}s=H_{0},$ and%
\begin{eqnarray*}
d^{\prime }\circ \left( \sigma \otimes \tau \right) &=&d^{\prime }\circ
\left( \alpha \left\langle sc,tc,F_{1}\right\rangle \cdot m^{\prime
}\left\langle 1_{sc},\tau \right\rangle \cdot \alpha ^{-1}\left\langle
sc,G_{1},td\right\rangle \cdot m^{\prime }\left\langle \sigma
,1_{td}\right\rangle \cdot \alpha \left\langle H_{1},sd,td\right\rangle
\right) \\
&=&\left( d^{\prime }\circ \alpha \left\langle sc,tc,F_{1}\right\rangle
\right) \cdot \left( d^{\prime }\circ m^{\prime }\left\langle 1_{sc},\tau
\right\rangle \right) \cdot \\
&&\ \ \ \ \ \ \ \ \ \ \ \ \ \ \ \ \ \ \ \ \left( d^{\prime }\circ \alpha
^{-1}\left\langle sc,G_{1},td\right\rangle \right) \cdot \left( d^{\prime
}\circ m^{\prime }\left\langle \sigma ,1_{td}\right\rangle \right) \cdot
\left( d^{\prime }\circ \alpha \left\langle H_{1},sd,td\right\rangle \right)
\\
&=&1_{d^{\prime }F_{1}}\cdot 1_{d^{\prime }F1}\cdot 1_{d^{\prime }td}\cdot
1_{d^{\prime }td}\cdot 1_{d^{\prime }td} \\
&=&1_{d^{\prime }F_{1}}\cdot 1_{d^{\prime }td}=1_{d^{\prime }F_{1}}\cdot
1_{F_{0}d}=1_{d^{\prime }F_{1}}\cdot 1_{d^{\prime }F_{1}}=1_{d^{\prime
}F_{1}}
\end{eqnarray*}%
with similar computations for $c^{\prime }\circ \left( \sigma \otimes \tau
\right) =1_{c^{\prime }H_{1}}$.

Commutativity of diagrams $\left( \ref{Octogon}\right) $ and $\left( \ref%
{pentagon tau}\right) $ is obtained using Yoneda Lemma, writing the
respective diagrams and adding all the possible arrows to fill them in order
to obtain the following \emph{mask} and \emph{diamond}%
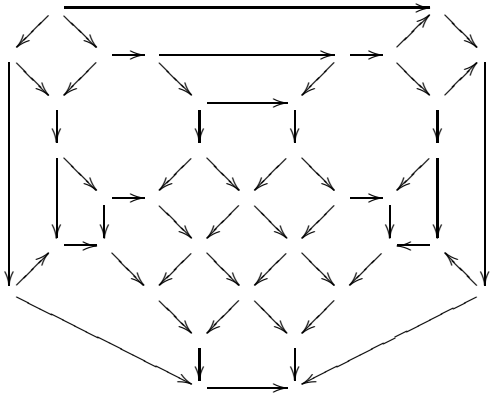
\begin{equation}
\ \ \ \ \ \ \ \ \ \ 
\xymatrix@=1pc{
& 
	\ar[rrrrrrrr]
	\ar[dr]
	\ar[dl]
& & & & & & & &
	\ar[rd]
&\\
	\ar[rd]
	\ar[ddddd]
&  & 
	\ar[r]
	\ar[ld]
&
	\ar[rrrr]
	\ar[rd]
& & & &
	\ar[r]
	\ar[ld]
&
	\ar[rd]
	\ar[ur]
&  &					
	\ar[ddddd]
& & &\\
&
	\ar[d]
& & &
	\ar[d]
	\ar[rr]
& &
	\ar[d]
& & &
	\ar[d]
	\ar[ru]
&&&\\
&
	\ar[rd]
	\ar[dd]
& & &
	\ar[dl]
	\ar[dr]
& &
	\ar[dl]
	\ar[dr]
& & &
	\ar[dl]
	\ar[dd]
\\
& &
	\ar[d]
	\ar[r]
&
	\ar[dr]
& &
	\ar[dl]
	\ar[dr]
& &
	\ar[dl]
	\ar[r]
&
	\ar[d]
\\
&
	\ar[r]
&
	\ar[rd]
 & &
	\ar[dl]
	\ar[dr]
& &
	\ar[dl]
	\ar[dr]
& &
	\ar[dl]
&	
	\ar[l]
\\
	\ar[rrrrdd]
	\ar[ru]
& & & 
	\ar[dr]
& &
	\ar[dl]
	\ar[dr]
& &
	\ar[dl]
& & &
	\ar[lllldd]
	\ar[lu]
&&&&&\\
& & & &
	\ar[d]
& &
	\ar[d]
\\
& & & &
	\ar[rr]
&   &
}
\tag{mask}
\end{equation}%
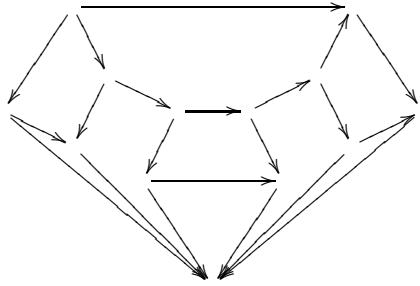
\begin{equation}
\xymatrix@=0.6pc{
& &
	\ar[rrrrrrrr]
	\ar[dddll]
	\ar[rdd]
& & & & & & & &
	\ar[rrddd]
&&\\
\\
& & &
	\ar[ddl]
	\ar[rrd]
& & & & & &
	\ar[rdd]
	\ar[ruu]
&&&\\
	\ar[rrd]
	\ar[dddddrrrrrr]
& & & & &
	\ar[rr]
	\ar[ddl]
& &
	\ar[rru]
	\ar[rdd]
& & & & &
	\ar[dddddllllll]
\\
& &
	\ar[ddddrrrr]
& & & & & & & &
	\ar[urr]
	\ar[ddddllll]
&&\\
& & & &
	\ar[rrrr]
	\ar[dddrr]
& & & &
	\ar[dddll]
&&&&\\
\\
\\
& & & & & &
&&&&&&\\
}
\tag{diamond}
\end{equation}%
in which squares commute by naturality, hexagons commute by definition of $%
\left( \sigma \otimes \tau \right) $, octagons commute because $S,T$ are
pseudo-natural transformations, pentagons in the diamond commute by the same
reason and all the other pentagons and triangles commute by coherence.
\end{proof}

\bigskip

The horizontal composition of pseudo-natural transformations is only defined
up to an isomorphism and it will be considered at the end of this paper.

In the next section\ we define square pseudo-modification ( simply called
pseudo-modification) and show that given two pseudo-categories $C,C^{\prime
} $, we obtain a pseudo-category by considering the pseudo-functors as
objects, natural transformations as morphisms, pseudo-natural
transformations as pseudo-morphisms and pseudo-modifications as cells. So,
in particular, we will show that the vertical composition of pseudo-natural
transformations is associative and has identities up to isomorphism. We also
show that PsCat is Cartesian closed up to isomorphism, that is, instead of
an isomorphism of categories PsCat($A\times B,C$)$\cong $PsCat($A,$PsCAT($%
B,C $)) we get an equivalence of categories PsCat($A\times B,C$)$\sim $PsCat(%
$A,$PsCAT($B,C$)).

\section{Pseudo-modifications}

Let \textbf{C} be a 2-category. Suppose $C,C^{\prime }$ are
pseudo-categories in \textbf{C}, $F,G,H,K:C\longrightarrow C^{\prime }$ are
pseudo-functors, $T=\left( t,\tau \right) :F\longrightarrow G$, $T^{\prime
}=\left( t^{\prime },\tau ^{\prime }\right) :H\longrightarrow K$ are
pseudo-natural transformations and $\theta =\left( \theta _{0},\theta
_{1}\right) :F\longrightarrow H,\theta ^{\prime }=\left( \theta _{0}^{\prime
},\theta _{1}^{\prime }\right) :G\longrightarrow K$ are two natural
transformations.

A \textbf{pseudo-modification} $\Phi $ (that will be represented as)%
\begin{equation*}
\squarecell{F}   {T}    {G}
                  {\theta}   {\Phi}    {\theta^\prime}
                  {H}   {T^\prime}    {K}%
\end{equation*}%
is a 2-cell of \textbf{C}%
\begin{equation*}
\Phi :t\longrightarrow t^{\prime }
\end{equation*}%
satisfying%
\begin{eqnarray}
d^{\prime }\circ \Phi &=&\theta _{0}  \label{dom_cod_ for phi} \\
c^{\prime }\circ \Phi &=&\theta _{0}^{\prime }  \notag
\end{eqnarray}%
and the commutativity of the square%
\begin{equation}
\begin{tabular}{ll}
$%
\begin{CD}%
\newline
\bullet 
@>%
\tau 
>%
>%
\bullet 
\\ %
\newline
@V%
m^{\prime }\left\langle \theta _{1}^{\prime },\Phi \circ d\right\rangle 
V%
V%
@V%
V%
m^{\prime }\left\langle \Phi \circ c,\theta _{1}\right\rangle 
V%
\\ %
\newline
\bullet 
@>%
\tau ^{\prime }%
>%
>%
\bullet 
&%
&%
\\ %
\newline
\end{CD}%
$ & \ \ \ \ \ \ .%
\end{tabular}
\label{square phi}
\end{equation}

Consider the case where \textbf{C}=Cat. Suppose $W,W^{\prime }$ are two
pseudo-categories in Cat, $F,G,H,K:W\longrightarrow W^{\prime }$ are
pseudo-functors, $T:F\longrightarrow G,T^{\prime }:H\longrightarrow K$ are
pseudo-natural transformations and $\theta :F\longrightarrow G,\theta
^{\prime }:H\longrightarrow K$ are natural transformations.

A pseudo-modification $\Phi $%
\begin{equation*}
\squarecell{F}   {T}    {G}
                  {\theta}   {\Phi}    {\theta^\prime}
                  {H}   {T^\prime}    {K}%
\end{equation*}%
is a family of cells%
\begin{equation*}
\squarecell{FA}   {t_A}    {GA}
                  {\theta_A}   {\Phi_A}    {\theta^{\prime}_{A}}
                  {HA}   {t^{\prime}_{A}}    {KA}%
\end{equation*}%
of $W^{\prime }$, for each object $A$ in $W$, where the square%
\begin{equation}
\begin{tabular}{ll}
$%
\begin{CD}%
\newline
t_{A}%
@>%
\Phi _{A}%
>%
>%
t_{A}^{\prime }%
\\ %
\newline
@V%
t_{a}%
V%
V%
@V%
V%
t_{a}^{\prime }%
V%
\\ %
\newline
t_{A^{\prime }}%
@>%
\Phi _{A^{\prime }}%
>%
>%
t_{A^{\prime }}^{\prime }%
&%
&%
\\ %
\newline
\end{CD}%
$ & ,%
\end{tabular}
\label{square tetra 1}
\end{equation}%
commutes for every morphism $a:A\longrightarrow A^{\prime }$ in $W$
(naturality of $\Phi $) and the square%
\begin{equation}
\begin{tabular}{ll}
$%
\begin{CD}%
\newline
Gf\otimes t_{A}%
@>%
\tau _{f}%
>%
>%
t_{B}\otimes Ff%
\\ %
\newline
@V%
\theta _{f}^{\prime }\otimes \Phi _{A}%
V%
V%
@V%
V%
\Phi _{B}\otimes \theta _{f}%
V%
\\ %
\newline
Kf\otimes t_{A}^{\prime }%
@>%
>%
\tau _{f}^{\prime }%
>%
t_{B}^{\prime }\otimes Hf%
&%
&%
\\ %
\newline
\end{CD}%
$ & ,%
\end{tabular}
\label{square tetra 2}
\end{equation}%
commutes for every pseudo-morphism $f:A\longrightarrow B$ in $W$.

Both squares $\left( \ref{square tetra 1}\right) $ and $\left( \ref{square
tetra 2}\right) $ may be displayed together with full information, for a $%
\varphi $ in $W$, as follows%
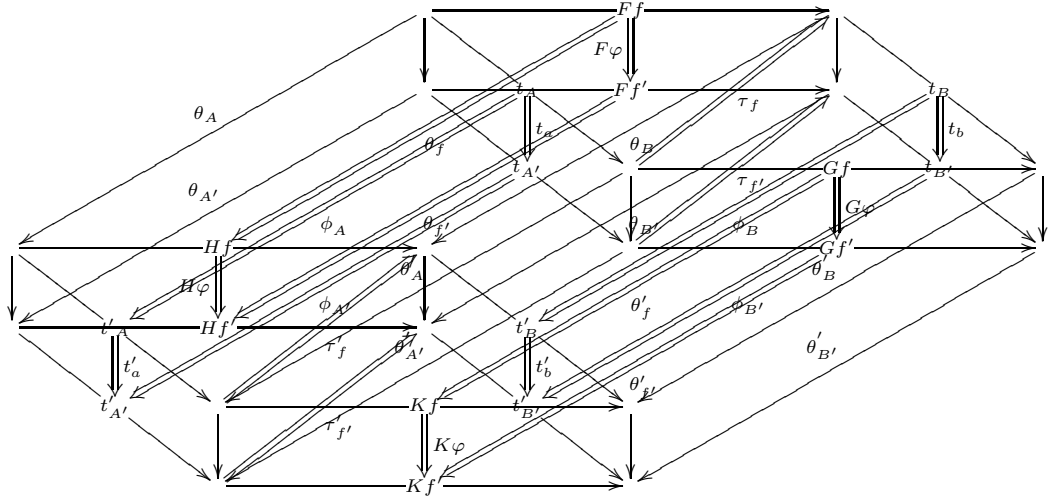
\begin{equation}
\xymatrix@C=6pc{
& &
	\ar[rr]|{Ff}="Ff"
	\ar[d]_{}
	\ar[dddll]_{\theta_A}
	\ar[rdd]|{t_A}="tA"
& &
	\ar[d]^{}
	\ar[dddll]^{\theta_B}
	\ar[rdd]|{t_B}="tB"
&\\
& &
	\ar[rr]|{Ff^{\prime}}="Ffp"
	\ar[dddll]_{\theta_{A^\prime}}
	\ar[rdd]|{t_{A^\prime}}="tAp"
& &
	\ar[rdd]|{t_{B^\prime}}="tBp"
	\ar[dddll]^{\theta_{B^\prime}}
&\\
& & &
	\ar[rr]|{Gf}="Gf"
	\ar[d]_{}
	\ar[dddll]_{\theta^{^\prime}_{A}}
	\ar@{=>}[uur]_{\tau_f}
& &
	\ar[d]^{}
	\ar[dddll]_{\theta^{^\prime}_{B}}
\\
	\ar[rr]|{Hf}="Hf"
	\ar[d]_{}
	\ar[rdd]|{t{^\prime}_A}="tpA"
& &
	\ar[d]^{}
	\ar[rdd]|{t^{\prime}_B}="tpB"
&
	\ar[rr]|{Gf^\prime}="Gfp"
	\ar[dddll]_{\theta^{^\prime}_{A^\prime}}
	\ar@{=>}[uur]_{\tau_{f^\prime}}
& & 
	\ar[dddll]_{\theta^{^\prime}_{B^\prime}}
\\
	\ar[rr]|{Hf^{\prime}}="Hfp"
	\ar[rdd]|{t^{\prime}_{A^\prime}}="tpAp"
& &
	\ar[rdd]|{t^{\prime}_{B^\prime}}="tpBp"
&&&\\
&
	\ar[rr]|{Kf}="Kf"
	\ar[d]_{}
	\ar@{=>}[uur]_{\tau^{\prime}_f}
& &
	\ar[d]^{}
&&\\
&
	\ar[rr]|{Kf^\prime}="Kfp"
	\ar@{=>}[uur]_{\tau^{\prime}_{f^\prime}}
& & & &\\
\ar@{=>}^{\theta_f} "Ff";"Hf" 
\ar@{=>}^{\theta_{f^\prime}} "Ffp";"Hfp" 
\ar@{=>}^{\theta^{\prime}_f} "Gf";"Kf" 
\ar@{=>}^{\theta^{\prime}_{f^\prime}} "Gfp";"Kfp" 
\ar@{=>}_{F\varphi} "Ff";"Ffp" 
\ar@{=>}^{G\varphi} "Gf";"Gfp" 
\ar@{=>}_{H\varphi} "Hf";"Hfp" 
\ar@{=>}^{K\varphi} "Kf";"Kfp"
\ar@{=>}^{t_a} "tA";"tAp"
\ar@{=>}^{t_b} "tB";"tBp"
\ar@{=>}^{t^{\prime}_a} "tpA";"tpAp"
\ar@{=>}^{t^{\prime}_b} "tpB";"tpBp"
\ar@{=>}^{\phi_A} "tA";"tpA"
\ar@{=>}^{\phi_{A^\prime}} "tAp";"tpAp"
\ar@{=>}^{\phi_B} "tB";"tpB"
\ar@{=>}^{\phi_{B^\prime}} "tBp";"tpBp"
}
\label{tetra picture}
\end{equation}

Return to the general case.

Let \textbf{C} be a 2-category and consider $C,C^{\prime }$ two
pseudo-categories in \textbf{C} as in $\left( \ref{C and C' again}\right) $.
Suppose $T,T^{\prime },T^{\prime \prime }$ are pseudo-natural
transformations between pseudo-functors from $C$ to $C^{\prime }$: we define
for%
\begin{equation*}
T\overset{\Phi }{\longrightarrow }T^{\prime }\overset{\Phi ^{\prime }}{%
\longrightarrow }T^{\prime \prime }
\end{equation*}%
a composition $\Phi ^{\prime }\Phi $ as the composition of 2-cells in 
\textbf{C}, and clearly it is well defined, is associative and has
identities. Now for $\theta ,\theta ^{\prime },\theta ^{\prime \prime }$
natural transformations between pseudo-functors from $C$ to $C^{\prime }$,
we define for 
\begin{equation*}
\theta \overset{\Phi }{\longrightarrow }\theta ^{\prime }\overset{\Psi }{%
\longrightarrow }\theta ^{\prime \prime }
\end{equation*}%
a pseudo-composition $\Psi \otimes \Phi =m^{\prime }\left\langle \Psi ,\Phi
\right\rangle .$

\begin{proposition}
\label{Prop well defined wcomp}Let \textbf{C} be a 2-category and suppose $%
\Psi ,\Phi $ are pseudo-modifications%
\begin{equation*}
\xymatrix@R=1pc@C=5pc{
F \ar[r]|{T}="f" \ar[d]_{\theta} & G \ar[d]^{\theta^\prime} \ar[r]|{S}="g" & H \ar[d]^{\theta^{\prime\prime}} \\
F^{\prime} \ar[r]|{T^\prime}="fp" & G^{\prime}  \ar[r]|{S^\prime}="gp" & H^{\prime}\\
\ar@{=>}^{\Phi} "f";"fp"
\ar@{=>}^{\Psi} "g";"gp"}%
\end{equation*}%
with $F,G,H,F^{\prime },G^{\prime },H^{\prime }$ pseudo-functors from $C$ to 
$C^{\prime }$ (pseudo-categories as in $\left( \ref{C and C' again}\right) $%
), $S,T,S^{\prime },T^{\prime }$ pseudo-natural transformations and $\theta
,\theta ^{\prime },\theta ^{\prime \prime }$ natural transformations as
considered above.\newline
The \ formula%
\begin{equation*}
\Psi \otimes \Phi =m^{\prime }\left\langle \Psi ,\Phi \right\rangle
\end{equation*}%
for pseudo-composition of pseudo-modifications is well defined.
\end{proposition}

\begin{proof}
Recall that the composition of pseudo-modifications is given by%
\begin{equation*}
S\otimes T=\left( m^{\prime }\left\langle s,t\right\rangle ,\left( \sigma
\otimes \tau \right) \right)
\end{equation*}%
with $\left( \sigma \otimes \tau \right) $ given as in $\left( \ref%
{tau_sigma}\right) $, hence%
\begin{equation*}
m^{\prime }\left\langle \Psi ,\Phi \right\rangle :m^{\prime }\left\langle
s,t\right\rangle \longrightarrow m^{\prime }\left\langle s^{\prime
},t^{\prime }\right\rangle
\end{equation*}%
is a 2-cell of \textbf{C} as required.\newline
Conditions $\left( \ref{dom_cod_ for phi}\right) $ are satisfied,%
\begin{eqnarray*}
d^{\prime }m^{\prime }\circ \left\langle \Psi ,\Phi \right\rangle
&=&d^{\prime }\pi _{2}^{\prime }\circ \left\langle \Psi ,\Phi \right\rangle
=d^{\prime }\circ \Phi =\theta _{0} \\
c^{\prime }m\circ \left\langle \Psi ,\Phi \right\rangle &=&c^{\prime }\pi
_{1}^{\prime }\left\langle \Psi ,\Phi \right\rangle =c^{\prime }\circ \Psi
=\theta _{0}^{\prime \prime }.
\end{eqnarray*}%
To prove commutativity of square $\left( \ref{square phi}\right) $ we use
Yoneda Lemma and the following diagram, obtained by adapting $\left( \ref%
{square phi}\right) $ to the present case and filling its interior%
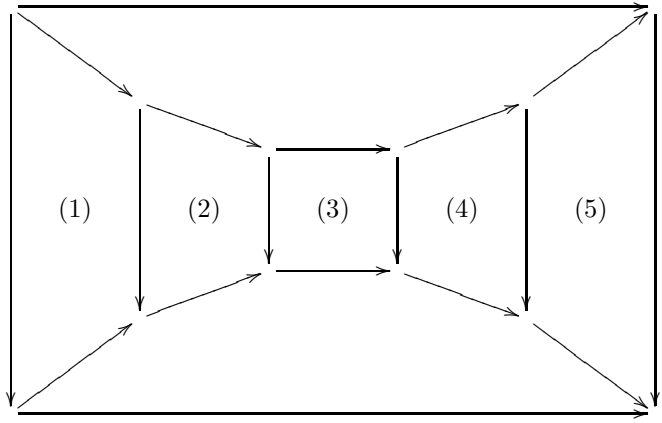
\begin{equation*}
\xymatrix@=1pc{
	\ar[rrrrrrrrrr]
	\ar[rrdd]
	\ar[dddddddd]
&&&&&&&&&&
	\ar[dddddddd]
\\
\\
&&
	\ar[rrd]
	\ar[dddd]
&&&&&&
	\ar[rruu]
	\ar[dddd]
&&\\
&&&&
	\ar[rr]
	\ar[dd]
&&
	\ar[rru]
	\ar[dd]
&&&&\\
& (1) && (2) && (3) && (4) && (5) & 
\\
&&&&
	\ar[rr]
&&
	\ar[rrd]
&&&\\
&&
	\ar[rru]
&&&&&&
	\ar[rrdd]
&&\\
\\
	\ar[rrrrrrrrrr]
	\ar[rruu]
&&&&&&&&&&\\
}%
\end{equation*}%
where hexagons commute by definition of $\left( \sigma \otimes \tau \right) $
and $\left( \sigma ^{\prime }\otimes \tau ^{\prime }\right) $, squares $%
\left( 1\right) ,\left( 3\right) ,\left( 5\right) $ commute by naturality of 
$\alpha ^{\prime }$ while squares $\left( 2\right) ,\left( 4\right) $
commute because $\Psi ,\Phi $ are pseudo-modifications (satisfy $\left( \ref%
{square tetra 2}\right) $) together with the fact that pseudo-composition
(in $C^{\prime }$) satisfies the middle interchange law $\left( \ref%
{interchange law}\right) $.
\end{proof}

Composition of pseudo-natural transformations is not associative, however
there is a special pseudo-modification for each triple of composable
pseudo-natural transformations.

\begin{proposition}
\label{Prop assoc}Let \textbf{C} be 2-category and suppose $%
F,G,H,K:C\longrightarrow C^{\prime }$ are pseudo-functors in \textbf{C} and
that $S=\left( s,\sigma \right) ,T=\left( t,\tau \right) ,U=\left(
u,\upsilon \right) $ are pseudo-natural transformations as follows%
\begin{equation*}
F\overset{S}{\longrightarrow }G\overset{T}{\longrightarrow }H\overset{U}{%
\longrightarrow }K.
\end{equation*}%
The 2-cell $\alpha _{U,T,S}^{\prime }=\alpha ^{\prime }\left\langle
u,t,s\right\rangle $ is a pseudo-modification \ 
\begin{equation*}
\squarecell{F}   {U\otimes (T\otimes S)}    {K}
                  {1}   {\alpha^{\prime}_{U,T,S}}    {1{\ \ ,}}
                  {F}   {(U\otimes T)\otimes S}    {K}%
\end{equation*}%
and it is natural in $S,T,U,$ in the sense that the square%
\begin{equation*}
\begin{tabular}{ll}
$%
\begin{CD}%
\newline
U\otimes \left( T\otimes S\right) 
@>%
\alpha ^{\prime }\left\langle u,t,s\right\rangle 
>%
>%
\left( U\otimes T\right) \otimes S%
\\ %
\newline
@V%
\varphi \otimes \left( \gamma \otimes \delta \right) 
V%
V%
@V%
V%
\left( \varphi \otimes \gamma \right) \otimes \delta 
V%
\\ %
\newline
U^{\prime }\otimes \left( T^{\prime }\otimes S^{\prime }\right) 
@>%
\alpha ^{\prime }\left\langle u^{\prime },t^{\prime },s^{\prime
}\right\rangle 
>%
>%
\left( U^{\prime }\otimes T^{\prime }\right) \otimes S^{\prime }%
&%
&%
\\ %
\newline
\end{CD}%
$ & 
\end{tabular}%
\end{equation*}%
commutes for every pseudo-modification $\varphi :U\longrightarrow U^{\prime
},\gamma :T\longrightarrow T^{\prime },\delta :S\longrightarrow S^{\prime }$.
\end{proposition}

\begin{proof}
The 2-cell $\alpha ^{\prime }\left\langle u,t,s\right\rangle $ is obtained
from%
\begin{equation*}
C_{0}\overset{\left\langle u,t,s\right\rangle }{\longrightarrow }%
C_{1}^{\prime }\times _{C_{0}^{\prime }}C_{1}^{\prime }\times
_{C_{0}^{\prime }}C_{1}^{\prime }\underrightarrow{\overrightarrow{\ \ \ \
\Downarrow \alpha ^{\prime }\ \ \ }}C_{1},
\end{equation*}%
and 
\begin{eqnarray*}
U\otimes \left( T\otimes S\right) &=&\left( m\left( 1\times m\right)
\left\langle u,t,s\right\rangle ,\left( \upsilon \otimes \left( \tau \otimes
\sigma \right) \right) \right) \\
\left( U\otimes T\right) \otimes S &=&\left( m\left( m\times 1\right)
\left\langle u,t,s\right\rangle ,\left( \left( \upsilon \otimes \tau \right)
\otimes \sigma \right) \right) ,
\end{eqnarray*}%
hence%
\begin{equation*}
\alpha ^{\prime }\left\langle u,t,s\right\rangle :m\left( 1\times m\right)
\left\langle u,t,s\right\rangle \longrightarrow m\left( m\times 1\right)
\left\langle u,t,s\right\rangle
\end{equation*}%
is a 2-cell of \textbf{C}.\newline
Conditions $\left( \ref{dom_cod_ for phi}\right) $ are satisfied%
\begin{eqnarray*}
d^{\prime }\circ \alpha ^{\prime }\circ \left\langle u,t,s\right\rangle
&=&1_{d^{\prime }\pi _{3}^{\prime }}\left\langle u,t,s\right\rangle
=1_{d^{\prime }s}=1_{F_{0}} \\
c^{\prime }\circ \alpha ^{\prime }\circ \left\langle u,t,s\right\rangle
&=&1_{c^{\prime }\pi _{1}^{\prime }}\left\langle u,t,s\right\rangle
=1_{c^{\prime }u}=1_{K_{0}}.
\end{eqnarray*}%
Commutativity of $\left( \ref{square phi}\right) $ follows from Yoneda Lemma
and the commutativity of the following diagram%
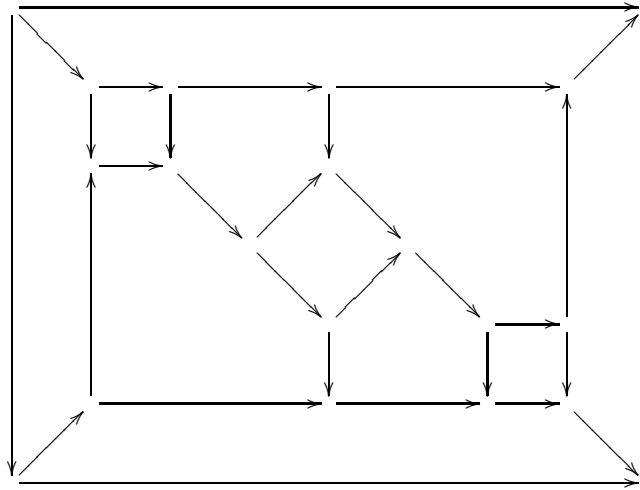
\begin{equation*}
\xymatrix{
	\ar[rrrrrrrr]
	\ar[dddddd]
	\ar[rd]
&&&&&&&&
	\ar[dddddd]
\\
&
	\ar[r]
	\ar[d]
&
	\ar[rr]
	\ar[d]
&&
	\ar[rrr]
	\ar[d]
&&&
	\ar[ru]
&\\
&
	\ar[r]
&
	\ar[rd]
&&
	\ar[rd]
&&\\
&&&
	\ar[rd]
	\ar[ru]
&&
	\ar[rd]
&&&\\
&&&&
	\ar[ru]
	\ar[d]
&&
	\ar[r]
	\ar[d]
&
	\ar[uuu]
	\ar[d]
&\\
&
	\ar[uuu]
	\ar[rrr]
&&&
	\ar[rr]
&&
	\ar[r]
&
	\ar[rd]
&\\
	\ar[ru]
	\ar[rrrrrrrr]
&&&&&&&&\\
}%
\end{equation*}%
where hexagons commute because $S,T,U$ are pseudo-natural transformations,
squares commute by naturality and pentagons by coherence.

To prove naturality we observe that%
\begin{eqnarray*}
\left( \left( \varphi \otimes \gamma \right) \otimes \delta \right) \cdot
\left( \alpha ^{\prime }\left\langle u,t,s\right\rangle \right) &=&\left(
m^{\prime }\left\langle m\left\langle \varphi ,\gamma \right\rangle ,\delta
\right\rangle \right) \cdot \left( \alpha ^{\prime }\left\langle
u,t,s\right\rangle \right) \\
&=&\left( m^{\prime }\left( m^{\prime }\times 1\right) \left\langle \varphi
,\gamma ,\delta \right\rangle \right) \cdot \left( \alpha ^{\prime
}\left\langle u,t,s\right\rangle \right) \\
&=&\left( 1_{m^{\prime }\left( m^{\prime }\times 1\right) }\cdot \alpha
^{\prime }\right) \circ \left( \left\langle \varphi ,\gamma ,\delta
\right\rangle \cdot 1_{\left\langle u,t,s\right\rangle }\right) \\
&=&\alpha ^{\prime }\circ \left\langle \varphi ,\gamma ,\delta \right\rangle
\end{eqnarray*}%
and%
\begin{eqnarray*}
\left( \alpha ^{\prime }\left\langle u^{\prime },t^{\prime },s^{\prime
}\right\rangle \right) \cdot \left( \varphi \otimes \left( \gamma \otimes
\delta \right) \right) &=&\left( \alpha ^{\prime }\left\langle u^{\prime
},t^{\prime },s^{\prime }\right\rangle \right) \cdot \left( m^{\prime
}\left\langle \varphi ,m^{\prime }\left\langle \gamma ,\delta \right\rangle
\right\rangle \right) \\
&=&\left( \alpha ^{\prime }\left\langle u^{\prime },t^{\prime },s^{\prime
}\right\rangle \right) \cdot \left( m^{\prime }\left( 1\times m^{\prime
}\right) \left\langle \varphi ,\gamma ,\delta \right\rangle \right) \\
&=&\left( \alpha ^{\prime }\cdot 1_{m^{\prime }\left( 1\times m^{\prime
}\right) }\right) \circ \left( 1_{\left\langle u^{\prime },t^{\prime
},s^{\prime }\right\rangle }\left\langle \varphi ,\gamma ,\delta
\right\rangle \right) \\
&=&\alpha ^{\prime }\circ \left\langle \varphi ,\gamma ,\delta \right\rangle
.
\end{eqnarray*}
\end{proof}

For every pseudo-functor there is a pseudo-identity pseudo-natural
transformation and a pseudo-identity pseudo-modification.

\begin{proposition}
\label{Prop idF}Consider a pseudo-functor $F=\left( F_{0},F_{1},\mu
_{F},\varepsilon _{F}\right) :C\longrightarrow C^{\prime }$ in a 2-category 
\textbf{C} (with $C,C^{\prime }$ pseudo-categories in \textbf{C} as in $%
\left( \ref{C and C' again}\right) $). The pair $\left( e^{\prime
}F_{0},\left( \lambda ^{\prime -1}\rho ^{\prime }\right) \circ F_{1}\right) $
is a pseudo-natural transformation in PsCat(\textbf{C})%
\begin{equation*}
id_{F}=\left( e%
{\acute{}}%
F_{0},\lambda ^{\prime -1}\rho ^{\prime }F_{1}\right) :F\longrightarrow F,
\end{equation*}%
and the 2-cell $1_{e^{\prime }F_{0}}:e^{\prime }F_{0}\longrightarrow
e^{\prime }F_{0}$ is a \ pseudo-modification in PsCat(\textbf{C})%
\begin{equation*}
\squarecell{F}   {id_F}    {F}
                  {1}   {1_{id_F}}    {1\ \ .}
                  {F}   {id_F}    {F}%
\end{equation*}
\end{proposition}

\begin{proof}
Clearly $e^{\prime }F_{0}:C_{0}\longrightarrow C_{1}^{\prime }$ is a
morphism of \textbf{C}, and 
\begin{equation*}
\lambda ^{\prime -1}\rho ^{\prime }F_{1}:m\left\langle F_{1},e^{\prime
}d^{\prime }F_{1}\right\rangle \longrightarrow m^{\prime }\left\langle
e^{\prime }c^{\prime }F_{1},F_{1}\right\rangle
\end{equation*}%
is a 2-cell (that is an isomorphism) of \textbf{C}.

Conditions $\left( \ref{dom cod t}\right) $ and $\left( \ref{dom cod tau}%
\right) $ are satisfied,%
\begin{eqnarray*}
d^{\prime }e^{\prime }F_{0} &=&F_{0} \\
c^{\prime }e^{\prime }F_{0} &=&F_{0}
\end{eqnarray*}%
\begin{eqnarray*}
d^{\prime }\circ \left( \lambda ^{\prime -1}\rho ^{\prime }F_{1}\right)
&=&d^{\prime }\circ \left( \lambda ^{\prime -1}\rho ^{\prime }\right) \circ
F_{1} \\
&=&\left( d^{\prime }\lambda ^{\prime -1}F_{1}\right) \cdot \left( d^{\prime
}\rho ^{\prime }F_{1}\right) \\
&=&\left( 1_{d^{\prime }F_{1}}\right) \cdot \left( 1_{d^{\prime
}F_{1}}\right) \\
&=&\left( 1_{d^{\prime }F_{1}}\right) ,
\end{eqnarray*}%
and similarly for $c^{\prime }\circ \left( \lambda ^{\prime -1}\rho ^{\prime
}F_{1}\right) =1_{c^{\prime }F_{1}}.$

Commutativity of $\left( \ref{Octogon}\right) $ is obtained using Yoneda
Lemma and the commutativity of the\ diagram%
\begin{equation*}
\xymatrix@C=0.0pc{
&
F_{g\otimes f} \otimes id_{FA}
	\ar[rr]
	\ar[rd]
	\ar[ld]
& &
id_{FC} \otimes F_{g\otimes f}
	\ar[rd]
&\\
(F_g\otimes F_f)\otimes id_{FA}
	\ar[rrd]
	\ar[d]
& &
F_{g\otimes f}
	\ar[ru]
	\ar[d]
& &
id_{FC} \otimes (F_g\otimes F_f)
	\ar[d]
\\
F_g\otimes (F_f\otimes id_{FA})
 	\ar[rd]
& &
F_g \otimes F_f
	\ar[ll]
	\ar[rru]
	\ar[rd]
& &
(id_{FC} \otimes F_g)\otimes F_f     \ \  ,
	\ar[ll]
\\
&
F_g\otimes ( id_{FB}\otimes F_f)
	\ar[ru]
	\ar[rr]
& &
(F_g\otimes  id_{FB})\otimes F_f
	\ar[ru]
&\\
}%
\end{equation*}%
while $\left( \ref{pentagon tau}\right) $ follows in a similar way as
observed in the diagram%
\begin{equation*}
\xymatrix{
&
F\left( id_{A}\right) \otimes id_{FA}
             \ar[rr]
             \ar[ld]
& &
id_{FA}\otimes F\left( id_{A}\right) 
             \ar[rd]
&\\
id_{FA}\otimes id_{FA}
             \ar[rrd]
& & & &
id_{FA}\otimes id_{FA}  .
             \ar[lld]
\\
& &
id_{FA}
}%
\end{equation*}

This proves that $id_{F}$ is a pseudo-natural transformation. To prove $%
1_{id_{F}}=1_{e^{\prime }F_{0}}$ is a pseudo-modification we note that%
\begin{equation*}
1_{e^{\prime }F_{0}}:e^{\prime }F_{0}\longrightarrow e^{\prime }F_{0}
\end{equation*}%
is a 2-cell of \textbf{C},%
\begin{eqnarray*}
d^{\prime }\circ 1_{e^{\prime }F_{0}} &=&1_{d^{\prime }e^{\prime
}F_{0}}=1_{F_{0}}, \\
c^{\prime }\circ 1_{e%
{\acute{}}%
F_{0}} &=&1_{c^{\prime }e^{\prime }F_{0}}=1_{F_{0}}.
\end{eqnarray*}%
To prove commutativity of square $\left( \ref{square phi}\right) $ we use
Yoneda Lemma and the commutativity of the following square%
\begin{equation*}
\begin{tabular}{ll}
$%
\begin{CD}%
\newline
Ff\otimes id_{FA}%
@>%
\lambda _{Ff}^{\prime -1}\rho _{Ff}^{\prime }%
>%
>%
id_{FB}\otimes Ff%
\\ %
\newline
@V%
1_{Ff}\otimes 1_{id_{FA}}%
V%
V%
@V%
V%
1_{id_{FB}}\otimes 1_{Ff}%
V%
\\ %
\newline
Ff\otimes id_{FA}%
@>%
\lambda _{Ff}^{\prime -1}\rho _{Ff}^{\prime }%
>%
>%
id_{FB}\otimes Ff%
&%
&%
\\ %
\newline
\end{CD}%
$ & .%
\end{tabular}%
\end{equation*}
\end{proof}

\begin{proposition}
\label{Prop lambda rho}Let \textbf{C} be a 2-category and suppose $%
F,G:C\longrightarrow C^{\prime }$ are pseudo-functors in \textbf{C}.\newline
For every pseudo-natural transformation%
\begin{equation*}
T=\left( t,\tau \right) :F\longrightarrow G
\end{equation*}%
there are two special pseudo-modifications%
\begin{equation*}
\squarecell{F}   { id_G \otimes T}    {G}
                  {1}   {\lambda_T}    {1   \ \ ,}
                  {F}   {T}    {G}%
\ \ 
\squarecell{F}   {T\otimes id_F}    {G}
                  {1}   {\rho_T}    {1   \ \ ,}
                  {F}   {T}    {G}%
\end{equation*}%
with $\lambda _{T}=\lambda ^{\prime }\circ t,\rho _{T}=\rho ^{\prime }\circ
t $ both natural in $T$.
\end{proposition}

\begin{proof}
It is clear that $\lambda ^{\prime }\circ t:m^{\prime }\left\langle
t,e^{\prime }F_{0}\right\rangle \longrightarrow t$ is a 2-cell of \textbf{C}%
, and%
\begin{eqnarray*}
d^{\prime }\circ \lambda ^{\prime }\circ t &=&1_{d^{\prime }t}=1_{F_{0}} \\
c^{\prime }\circ \lambda ^{\prime }\circ t &=&1_{c^{\prime }t}=1_{G_{0}}.
\end{eqnarray*}%
The commutativity of square $\left( \ref{square phi}\right) $ is obtained
from the commutativity of diagram%
\begin{equation*}
\xymatrix@C=-1pc{
G_f\otimes (id_{G_A}\otimes t_A)
	\ar[rrrrrr]
	\ar[rd]
	\ar[ddd]
& & & & & &
(id_{G_B}\otimes t_B)\otimes F_f
	\ar[ddd]
\\
&
(G_f\otimes id_{G_A})\otimes t_A
	\ar[rd]
	\ar[ldd]
& & & &
id_{G_B}\otimes (t_B\otimes F_f)
	\ar[ru]
	\ar[rdd]
&\\
& & 
(id_{G_B}\otimes G_f)\otimes t_A   - \ \ \ \ \ 
	\ar[rr]
	\ar[lld]
& &
id_{G_B}\otimes (G_f\otimes t_A)
	\ar[ru]
	\ar[lllld]
&&\\
G_f\otimes t_A
	\ar[rrrrrr]
&&&&&&
t_B\otimes F_f
\\
}%
\end{equation*}%
In order to prove naturality of $\lambda _{T}$ consider a internal
pseudo-modification 
\begin{equation*}
\squarecell{F}   {T}    {G}
                  {\theta}   {\Phi}    {\theta^\prime}
                  {H}   {T^\prime}    {K}%
\end{equation*}%
as defined in $\left( \ref{dom_cod_ for phi}\right) $; then, on the one hand
we have%
\begin{eqnarray*}
\Phi \cdot \left( \lambda ^{\prime }\circ t\right) &=&\left(
1_{C_{1}^{\prime }}\circ \Phi \right) \cdot \left( \lambda ^{\prime }\circ
1_{t}\right) \\
&=&\left( 1_{C_{1}^{\prime }}\cdot \lambda ^{\prime }\right) \circ \left(
\Phi \cdot 1_{t}\right) \\
&=&\lambda ^{\prime }\circ \Phi
\end{eqnarray*}%
and on the other hand we have%
\begin{eqnarray*}
\left( \lambda ^{\prime }\circ t^{\prime }\right) \cdot \left( m^{\prime
}\left\langle e^{\prime }\theta _{0}^{\prime },\Phi \right\rangle \right)
&=&\left( \lambda ^{\prime }\circ t^{\prime }\right) \cdot \left( m^{\prime
}\left\langle e^{\prime }c^{\prime },1_{C_{1}^{\prime }}\right\rangle \circ
\Phi \right) \\
&=&\left( \lambda ^{\prime }\cdot 1_{m^{\prime }\left\langle e^{\prime
}c^{\prime },1_{C_{1}^{\prime }}\right\rangle }\right) \circ \left(
1_{t^{\prime }}\cdot \Phi \right) \\
&=&\lambda ^{\prime }\circ \Phi .
\end{eqnarray*}

The proof on rho is similar.
\end{proof}

The three last propositions lead us to the following theorem.

\begin{theorem}
Let \textbf{C} be a 2-category, and consider $C,C^{\prime }$ two
pseudo-categories in \textbf{C}. The data:

\begin{itemize}
\item objects: pseudo-functors from $C$ to $C^{\prime }$;

\item morphisms: natural transformations (between pseudo-functors from $C$
to $C^{\prime }$);

\item pseudo-morphisms: pseudo-natural transformations (between
pseudo-functors from $C$ to $C^{\prime }$);

\item cells: pseudo-modifications (between such natural and pseudo-natural
transformations);
\end{itemize}

form a pseudo-category (in Cat).
\end{theorem}

\begin{proof}
Natural transformations and pseudo-functors form a category: theorem \ref%
{Theorem WCat is a 2Cat}. pseudo-modifications and pseudo-natural
transformations also form a category: the composition is associative and has
identities (that inherit the structure of 2-cells of the ambient 2-category).

For every pseudo-natural transformation $T=\left( t,\tau \right)
:F\longrightarrow G$, the identity pseudo-modification is $1_{T}=1_{t}$%
\begin{equation*}
\squarecell{F}   {T}    {G}
                  {1}   {1_T}    {1      \ \ .}
                  {F}   {T}    {G}%
\end{equation*}%
For each pair of pseudo-composable pseudo-modifications $\Phi ,\Psi ,$ there
is a (well defined - proposition \ref{Prop well defined wcomp})
pseudo-composition $\Phi \otimes \Psi =m^{\prime }\left\langle \Phi ,\Psi
\right\rangle $ satisfying $\left( \ref{interchange law}\right) $%
\begin{equation*}
\left( \Phi \Phi ^{\prime }\right) \otimes \left( \Psi \Psi ^{\prime
}\right) =m^{\prime }\left\langle \Phi \Phi ^{\prime },\Psi \Psi ^{\prime
}\right\rangle
\end{equation*}%
\begin{eqnarray*}
\left( \Phi \otimes \Psi \right) \left( \Phi ^{\prime }\otimes \Psi ^{\prime
}\right) &=&\left( m^{\prime }\left\langle \Phi ,\Psi \right\rangle \right)
\left( m^{\prime }\left\langle \Phi ^{\prime },\Psi ^{\prime }\right\rangle
\right) \\
&=&\left( 1_{m^{\prime }}1_{m^{\prime }}\right) \circ \left( \left\langle
\Phi ,\Psi \right\rangle \left\langle \Phi ^{\prime },\Psi ^{\prime
}\right\rangle \right) \\
&=&m^{\prime }\left\langle \Phi \Phi ^{\prime },\Psi \Psi ^{\prime
}\right\rangle ;
\end{eqnarray*}%
and $1_{T\otimes S}=1_{T}\otimes 1_{S}$,%
\begin{equation*}
1_{m\left\langle t,s\right\rangle }=1_{m}\circ 1_{\left\langle
t,s\right\rangle }=1_{m}\circ \left\langle 1_{t},1_{s}\right\rangle
=m\left\langle 1_{t},1_{s}\right\rangle .
\end{equation*}

For each natural transformation $\theta :F\longrightarrow G$ there is a
pseudo-modification%
\begin{equation*}
\squarecell{F}   {id_F}    {F}
                  {\theta}   {id_\theta}    {\theta}
                  {G}   {id_G}    {G}%
\end{equation*}%
with $id_{\theta }=e^{\prime }\theta _{0}$, satisfying%
\begin{equation*}
id_{1_{F}}=e^{\prime }1_{F_{0}}=1_{e^{\prime }F_{0}}=1_{id_{F}},
\end{equation*}%
\begin{equation*}
id_{\theta ^{\prime }\theta }=e^{\prime }\circ \left( \theta _{0}^{\prime
}\theta _{0}\right) =\left( e^{\prime }\circ \theta _{0}^{\prime }\right)
\left( e^{\prime }\circ \theta _{0}\right) =id_{\theta ^{\prime }}id_{\theta
}.
\end{equation*}%
By Proposition \ref{Prop assoc} there is a special pseudo-modification $%
\alpha _{T,U,S}=\alpha \left\langle T,U,S\right\rangle $ for each triple of
composable pseudo-natural transformations $T,U,S$, natural in each component
and satisfying the pentagon coherence condition.

By Proposition $\ref{Prop lambda rho}$ there are two special
pseudo-modifications $\lambda _{T},\rho _{T}$ to each pseudo-natural
transformation $T:F\longrightarrow G,$ natural in $T$ and satisfying the
triangle coherence condition.
\end{proof}

\section{Conclusion and final remarks}

The mathematical object PsCat that we have just defined has the following
structure:

\begin{itemize}
\item objects: $A,B,C,...$

\item morphisms: $f:A\longrightarrow B,...$

\item 2-cells: $\theta :f\longrightarrow g,...$($f,g:A\longrightarrow B$)

\item pseudo-cells: $%
\xymatrix{f \ar[r]|{T} & g}%
,...$

\item tetra cells: $%
\squarecell{f}   {T}    {g}
   {\theta}   {\Phi}    {\theta^\prime}
   {f^\prime}   {T^\prime}    {g^\prime}%
\begin{array}[t]{c}
\\ 
,...%
\end{array}%
$
\end{itemize}

where objects, morphisms and 2-cells form a 2-category and for each pair of
objects $A,B$, the morphisms, 2-cells, pseudo-cells and tetra cells from $A$
to $B$ form a pseudo-category.

Two questions arise at this moment:

- What is happening from PsCat$(B,C)\times $PsCat$(A,B)$ to PsCat$(A,C)$?

- What is the relation between PsCat$(A\times B,C)$ and PsCat$(A,$PsCAT$%
(B,C))$?

The answer to the second question is easy to find out. If starting with a
pseudo-functor in PsCat$\left( A\times B,C\right) $, say%
\begin{equation*}
h:A\times B\longrightarrow C,
\end{equation*}%
by going to PsCat$\left( A,C^{B}\right) $ and coming back we will obtain
either%
\begin{equation*}
h\left( c,g\right) \otimes h\left( f,b\right)
\end{equation*}%
or 
\begin{equation*}
h\left( f,d\right) \otimes h\left( a,g\right)
\end{equation*}%
instead of $h\left( f,g\right) $ as displayed in the diagram below%
\begin{equation*}
\xymatrix{
(a,b) 
	\ar[dd]^{(f,g)}
&  h(a,b)
	\ar[rr]^{h(f,b)}
	\ar[dd]_{h(a,g)}
	\ar[rrdd]^{h(f,g)}="h"
	\ar@{}[rrdd]_{}="h_"
& &  h(c,b)
	\ar[dd]^{h(c,g)}
	\ar@{ }[]="hcb"
\\
& &
\ar@{=>}[d]^{\mu}
\ar@{=>}[r]^{\mu}
&\\
(c,d) & h(a,d)
	\ar@{ }[]="had"
	\ar[rr]_{h(f,d)}
& & h(c,d)
\\
}%
\end{equation*}%
And since they are all isomorphic via $\mu $ and $\tau $ we have that the
relation is an equivalence of categories.

A similar phenomena happens when trying to define horizontal composition of
pseudo-natural transformations (while trying to answer the first question):
there are two equally good ways to define a horizontal composition and they
differ by an isomorphism.

Let \textbf{C} be a 2-category and $C,C^{\prime },C^{\prime \prime }$%
pseudo-categories in \textbf{C}, consider $S,T$ pseudo-natural
transformations as in%
\begin{equation*}
C\underset{G}{\overset{F}{\overrightarrow{\underrightarrow{\ \ \downarrow T\
\ }}}}C^{\prime }\underset{G^{\prime }}{\overset{F^{\prime }}{%
\overrightarrow{\underrightarrow{\ \ \downarrow S\ \ }}}}C^{\prime \prime }
\end{equation*}%
there are two possibilities to define horizontal composition%
\begin{equation*}
S\circ _{w1}T=m^{\prime \prime }\left\langle sG_{0},F_{1}t^{\prime
}\right\rangle
\end{equation*}%
and 
\begin{equation*}
S\circ _{w2}T=m^{\prime \prime }\left\langle G_{1}^{\prime
}t,sF_{0}\right\rangle
\end{equation*}%
as displayed in the following picture%
\begin{equation*}
\begin{tabular}{ll}
$%
\begin{CD}%
\newline
C_{1}%
@<%
e%
<%
<%
C_{0}%
\\ %
\newline
@V%
F_{1}%
V%
G_{1}%
V%
\hspace{-1cm}\swarrow _{t}\ \ \ 
@V%
F_{0}%
V%
G_{0}%
V%
\\ %
\newline
C_{1}^{\prime }%
@<%
e^{\prime }%
<%
<%
C_{0}^{\prime }%
\\ %
$\newline
$%
@V%
F_{1}^{\prime }%
V%
G_{1}^{\prime }%
V%
\hspace{-1cm}\swarrow _{s}\ \ \ 
@V%
F_{0}^{\prime }%
V%
G_{0}^{\prime }%
V%
\\ %
\newline
C_{1}^{\prime \prime }%
@<%
e^{\prime \prime }%
<%
<%
C_{0}^{\prime \prime }%
\\ %
\newline
\end{CD}%
$ & .%
\end{tabular}%
\end{equation*}

Hence we have two isomorphic functors from PsCat$(B,C)\times $PsCat$(A,B)$
to PsCat$(A,C)$ both defining a horizontal composition.

We note that this behaviour, of composition beeing defined up to
isomorphism, also occurs while trying to compose homotopies. So one can
expect further relations between the \ theory of pseudo-categories and
homotopy theory to be investigated.

For instance the category Top itself may be viewed as a \ structure with
objects (spaces), morphisms (continuous mappings), 2-cells (homotopy classes
of homotopies), pseudo-cells (simple homotopies) and tetra cells (homotopies
between homotopies).


\begin{thebibliography}{99}
\bibitem{Benabou} J. B\'{e}nabou, Introduction to Bicategories, Lecture
Notes in Mathematics, number 47, pages 1-77, Springer-Verlag, \ Berlin, 1967.

\bibitem{Gray} J.W. Gray, \textit{Formal Category Teory: Adjointness for
2-categories}, Lecture Notes in Mathematics, Springer-Verlag, 1974

\bibitem{Leinster} T. Leinster, \textit{Higher Operads, Higher Categories},
London Mathematical Society Lecture Notes Series, Cambridge University
Press, 2003 (electronic version).

\bibitem{ML} S. MacLane, \textit{Categories for the working Mathematician},
Springer-Verlag, 1998, 2$^{nd}$ edition.

\bibitem{MF1} N. Martins-Ferreira, Internal Bicategories in Ab, Preprint
CM03/I-24, Aveiro Universitiy, 2003.

\bibitem{MF2} N. Martins-Ferreira, Internal Weak Categories in Additive
2-Categories with Kernels, Fields Institute Communications, Volume 43,
p.387-410, 2004.

\bibitem{MF3} N. Martins-Ferreira, Weak categories in Grp, Unpublished.

\bibitem{Pare} R. Par\'{e} and M. Grandis, Adjoints for Double Categories,
Cahiers Topologie
et G\'{e}om\'{e}trie Difer\'{e}ntielle Cat\'{e}goriques, XLV(3),2004, 193-240.

\bibitem{Pare2} R. Par\'{e} and M. Grandis, Limits in Double Categories,
Cahiers Topologie
et G\'{e}om\'{e}trie Difer\'{e}ntielle Cat\'{e}goriques, XL(3), 1999, 162-220.

\bibitem{Power} J. Power, 2-Categories, BRICS Notes Series, NS-98-7.

\bibitem{Street} R.H. Street, Cosmoi of internal categories, Trans. Amer.
Math. Soc. 258, 1980, 271-318

\bibitem{Street2} R.H. Street, Fibrations in Bicategories, Cahiers Topologie
et G\'{e}om\'{e}trie Difer\'{e}ntielle Cat\'{e}goriques, 21:111-120, 1980.
\end{thebibliography}
\end{document}